\documentclass[aop,preprint,reqno]{imsart}
\RequirePackage{amsthm,amsmath,amsfonts,amssymb}
\RequirePackage[authoryear]{natbib}
\RequirePackage[colorlinks,citecolor=blue,urlcolor=blue]{hyperref}
\usepackage{enumitem}
\RequirePackage{graphicx,ulem}
\startlocaldefs
\theoremstyle{plain}

\newtheorem{theorem}{Theorem}[section]
\newtheorem{lemma}[theorem]{Lemma}
\newtheorem{proposition}[theorem]{proposition}
\newtheorem{corollary}[theorem]{Corollary}
\newtheorem{remark}{Remark}[section]
\theoremstyle{remark}
\newtheorem{Assumption}{Assumption}

\endlocaldefs

\begin{document}
		
\begin{frontmatter}
		\title{Revisit of a Diaconis urn model}
				
\begin{aug}
			\author{\fnms{Li} \snm{Yang}\ead[label=e1,mark]
				{yangl245@nenu.edu.cn}},
			\author{\fnms{Jiang} \snm{Hu}\ead[label=e2,mark]{huj156@nenu.edu.cn}},
			\and
			\author{\fnms{Zhidong} \snm{Bai}\ead[label=e3,mark]{baizd@nenu.edu.cn}},
			\address{KLASMOE and School of Mathematics and Statistics, Northeast Normal University, China.
				\printead{e1,e2,e3}}
		\end{aug}
				
\begin{abstract}
			Let $G$ be a finite Abelian group of order $d$. We consider an urn in which, initially, there are labeled balls that generate the group $G$. Choosing two balls from the urn with replacement, observe their labels, and perform a group multiplication on the respective group elements to obtain a group element. Then, we put a ball labeled with that resulting element into the urn. This model was formulated by P. Diaconis while studying a group theoretic algorithm called MeatAxe (\cite{r33}). \cite{r27} partially investigated this model. In this paper, we further investigate and generalize this model. More specifically, we allow a random number of balls to be drawn from the urn at each stage in the Diaconis urn model. For such a case, we verify that the normalized urn composition converges almost surely to the uniform distribution on the group $G$. Moreover, we obtain the asymptotic joint distribution of the urn composition by using the martingale central limit theorem.
		\end{abstract}
					
				
			
\end{frontmatter}
		
\section{Introduction}
	The urn model originating from \cite{r1} can be described as follows. Initially, the urn contains $N$ balls, of which there are $n_1$ white balls and $n_2$ black balls. At each stage, a ball is randomly drawn from the urn and returned along with $r$ other balls of the same color. Repeat the above process. This model is a famous probability model, and it is frequently used to simulate disease infection. In recent years, some deformations of this model have been extensively studied, which has led to the widespread use of the urn model ensemble in economics, informatics, biology, machine learning, adaptive designs in clinical trials, and other fields. We refer the reader to \cite{r2}, \cite{r3}, \cite{r4} and \cite{r5} for wider applications.
	\subsection{Classical P\'{o}lya urns}
	In classical urn models, one ball is drawn with replacement at each stage. Due to the sequential nature of the P\'{o}lya urn, many new urn models have been derived by redefining the rules of drawing and adding balls, making urn models widely applicable in clinical trials. A typical model is the generalized P\'{o}lya urn (GPU) model, also called the GFU (generalized Friedman's urn) in the literature, which can be described as follows. Initially, the urn contains $\mathbf{Y}_{0}=(Y_{01},\dots, Y_{0K})$ balls, where $Y_{0k}$ denotes the number of balls of type $k\in\{1,\dots, K\}$. In fact, each type of ball can be imagined as a treatment, and patients arrive in sequence. At stage $i$, a ball is drawn with replacement. If the drawn ball is of type $q\in\{1,\dots, K\}$, then assign the $q$th treatment to the $i$th patient. We then obtain a covariate $\xi_i$, which is associated with the response of the $i$th patient. In the meantime, we add $D_{qk}(i)$ balls of type $k\in\{1,\dots, K\}$ to the urn, where $D_{qk}(i)$ is some function of $\xi_i$. Repeat the above process $n$ stages, we then obtain the urn composition $\mathbf{Y}_{n}=(Y_{n1},Y_{n2},\dots, Y_{nK})$ and the patient assignment $\mathbf{N}_{n}=(N_{n1},N_{n2},\dots, N_{nK})$, where $N_{nk}$ represents the number of patients assigned to treatment $k$ in the first $n$ stages.
		
	For the GPU model, there are two concerns, the limit properties of $\mathbf{Y}_{n}$ and $\mathbf{N}_{n}$, they satisfy the following recursion:
	$$\mathbf{Y}_{n}=\mathbf{Y}_{n-1}+\mathbf{X}_{n}\mathbf{D}_{n},\ \ \ \ \ \mathbf{N}_{n}=\sum\nolimits_{i=1}^{n}\mathbf{X}_{i}. $$
	where $\mathbf{D}(i)=(D_{qk}(i))_{K\times K}$, $\mathbf{X}_{{i}}$ represents the result of the drawing, i.e., $\mathbf{X}_{i}=\mathbf{e}_{k},\ k=1,2,\dots, K$, if the ball drawn for the $i$th stage is of type $k$, and $\mathbf{e}_{k} $ is a $K$-dimensional zero vector except the $k$-th element, which is one. We call $\mathbf{D}(i)$ the adding rule and let $\mathbf{H}_{i}=\mathbb{E}(\mathbf{D}(i)|{\mathcal F}_{i-1})$ be the generating matrix, where ${\mathcal F}_{n}=\sigma\{\mathbf{Y}_{0},\{\mathbf{X}_{i}\}_{i=1}^{n}, and \{\mathbf{D}_{i}\}_{i=1}^{n}\}$. If $\mathbf{H}_{i}=\mathbf{H}$ for all $i\geq1$, then this model is said to be {\it homogeneous}; otherwise, $\{\mathbf{H}_{i}\}$ is assumed to have a limit $\mathbf{H}$, and $\mathbf{H}$ is often considered irreducible. 

	Actually, the asymptotic properties of $\mathbf{Y}_{n}$ and $\mathbf{N}_{n}$ are closely related to the eigenvalues and the left eigenvector corresponding to the maximum eigenvalue of $\mathbf{H}$. \cite{r6,r7} first considered the asymptotic properties of $\mathbf{Y}_{n}$ with a homogeneous generating matrix $\mathbf{H}$ by using the theory of continuous-time branching processes, and they conjectured that the normalized $\mathbf{N}_{n}$ would be asymptotically normal. Decades later, \cite{r8}, \cite{r9,r10}, \cite{r11}, \cite{r12} and \cite{r39} resolved the conjecture by martingales theory and Gaussian approximation methods with a nonhomogeneous generating matrix. For the asymptotic properties of the urn model with an unbalanced update or a nonlinear drawing rule, we refer to \cite{r24}, \cite{r25} and \cite{r26}. In addition, there are some other commonly useful models, such as the drop-the-loser (DL) model \citep{r13}, randomly reinforced urn
	model \citep{r13,r14,r15,r16}, and some particular $homogeneous$ urn models, for which their replacement matrices are reducible, but some of the diagonal subblocks are irreducible; see \cite{r36} and \cite{r37} for details.
		
\subsection{Multiple-drawing urns}
	For the case that multiple balls are drawn at each stage, the first general results for two-color urns were obtained by \cite{r19}, under the assumptions of {\it balance} (i.e., the number of balls added to the urn at each stage is a fixed integer), {\it affinity} (i.e., the conditional expectation of the number of white balls in the urn is linear) and {\it tenability} (i.e., the process of drawing and replacing balls is continued indefinitely). 

	For multicolor urns, the commonly used method is stochastic approximation (SA). The normalized urn composition can be written within a Robbins-Monro (RM) scheme, which is an iterative SA. 
	Partial results of the first-order and second-order convergences of the SA algorithm were introduced in detail in \cite{r20}, \cite{r32} and \cite{r21}. \cite{r31} established the link between SA theory and the urn model proposed by \cite{r9,r10} and \cite{r35}. \cite{r22} established some central limit theorems (CLTs) for the SA algorithm under the Lindeberg condition. Later, the results obtained by \cite{r19} were generalized in \cite{r23}, removing the {\it affinity} assumption. The model they studied can be described as follows: a fixed number of balls, $m$, are drawn at each stage, and the replacement rule depends on the color of the $m$ balls. The model is assumed to be balanced and tenable. After $n$ stages, $\mathbf{Z_{n}}$, which denotes the normalized urn composition, satisfies a recurrence formula:
	$$\mathbf{Z}_{n+1}=\mathbf{Z}_{n}+\frac{1}{T_{n+1}}(\mathbf{h}(\mathbf{Z}_{n})+\Delta\mathbf{M}_{n+1}+{\mathbf{r}}_{n+1}),$$
	where $\mathbf{h}:{\mathbb{R}}^{d}\to{\mathbb{R}}^{d}$ is a real vector-valued function related to the adding rule, $\{\Delta\mathbf{M}_{n},{\mathcal{F}}_{n};n\ge1\}$ is a martingale difference sequence, and $\{{\mathbf{r}}_{n}\}$ is a remaining item sequence. By using the results of \cite{r20} and \cite{r22} and by verifying the conditions of the convergence theorem for SA, the authors obtained the limiting distribution of the normalized urn composition. Notably, the almost sure (a.s.) convergence of $\mathbf{Z}_{n}$ needs to satisfy the following necessary conditions: $\mathbf{h}$ is uniformly bounded and admits a zero $\mathbf{Z}$ of $\mathbf{h}$ such that for a positive integer $N$,
	\begin{equation}  \label{112}
		\langle \mathbf{h}(\mathbf{Z}_{n}),\ \mathbf{Z}_{n}-\mathbf{Z}\rangle<0,\ \mbox{for \ all} \ n>N, 
	\end{equation}
	where we denote the scalar product of the vectors $\mathbf{h}(\mathbf{Z}_{n})$ and $\mathbf{Z}_{n}-\mathbf{Z}$ by $\langle \mathbf{h}(\mathbf{Z}_{n}),\ \mathbf{Z}_{n}-\mathbf{Z}\rangle$.
	If $\mathbf{Z}_{n}\to\mathbf{Z}$ a.s.,
	then the second-order convergence and the rate of convergence of $\mathbf{Z}_{n}$ depend on ${\rm Re}(\Lambda)$, where $\Lambda$ is the eigenvalue of $-\nabla\mathbf{h}(\mathbf{Z})$ with the smallest real part. If $\mathbf{h}$ has no explicit expression, the above conditions are difficult to verify. Most recently,
	\cite{r38} considered a two-color urn model with a random multiple-drawing and a random time-dependent addition matrix, and they proved asymptotic properties for the proportion of balls of a given color.
		
\subsection{Urn model proposed by P. Diaconis}
	Let $G$ be a finite group with $d$ elements. Initially, the urn contains some balls, which are labeled by the generators of the group. At each stage, two balls are drawn randomly from the urn with replacement, observe their labels, and perform group multiplication on these two elements. Then, we obtain an element of group $G$, a ball labeled with the element then added to the urn. Repeat the above stages $n$ stages, we are interested in the proportion of balls labeled by various elements in the urn when $n$ tends to infinity.		
	This model was formulated by P. Diaconis when studying a group theoretic algorithm called MeatAxe; we refer to \cite{r33} for more details. Then, P. Diaconis proposed this model in person to Z. D. Bai, one of the authors of this paper, many years ago.
			
	Although this is a multiple-drawing urn model, when $d$ is large it is difficult to verify whether the model satisfies the assumptions of the SA algorithm convergence given in \cite{r23}, due to the particularity of the group multiplication. Furthermore, the limit properties of the urn composition are related to whether the urn contains balls labeled as generators of group $G$ initially, which was not considered in \cite{r23}.
		
	Let $G^{*}$ be the simplex of probability distributions over $G$. Let $T:G^{*}\to G^{*}$ be a map of the simplex into itself. For this Diaconis urn model, when $d=2$, \cite{r27} proved that the normalized urn composition converges to a uniform distribution over group $G$; when $d>2$, they established the existence of a limit. However, since the map $T$ has multiple fixed points when $d>2$, their result does not imply that the limit is the uniform distribution over group $G$. Moreover, they did not obtain the growth rate of the urn compositions; thus, \cite{r27} did not completely solve the problem proposed by Diaconis. Later, \cite{r28} treated the Diaconis urn model as a random walk and found that if the urn contains balls whose labels generate the finite group $G$, then with probability 1, the limiting fraction of balls labeled by each group element approaches $1/|G|$, where $|G|$ denotes the order of group $G$.
		
	In this paper, we focus on the case in which the group $G$ in the Diaconis urn model is an Abelian group, and considered a more general case where a random number of balls $M_i$ are drawn at stage $i$, instead of drawing two balls, where $\{M_i,\ i=1,2,\dots\}$ are independently and identically distributed (i.i.d). 
	First, we consider the group $G=({\mathbb Z}_{d},+)$ and the addition group of residue classes modulo $d$. By constructing a subsequence $\{\mathbf{a}_{n_{k}},k=1,2,\dots\}$ of normalized urn composition $\{\mathbf{a}_{n},n=1,2,\dots\}$, we demonstrate that individual compositions of the subsequence converge to $1/|G|$ almost surely. Then, we obtain the result that when $k$ is large enough, for all $n_{k}<n<n_{k+1}$ and $j\in\{0,1,\dots, d-1\}$, the distance between $a_{nj}$ and $a_{n_{k},j}$ is small enough. Thus, we find that the individual compositions of $\{\mathbf{a}_{n}\}$ converge to $1/|G|$ almost surely when $n$ tends to infinity.
	Furthermore, we extend the convergence results for $({\mathbb Z}_{d},+)$ to those of arbitrary Abelian groups. Second, by applying the martingale CLT, we obtain the CLT for the normalized urn composition.	
	To the best of our knowledge, this is the first paper to establish the CLT of the Diaconis urn model, verifying Diaconis' conjecture about the second-order asymptotic properties of the normalized urn composition. And, we considered the results of drawing a random number of balls at each stage, the gap in the research of multicolor urn models with random multiple-drawing has been filled.
		
	The paper is organized as follows. In Section 2, we introduce the needed notation and state the main results. Proofs of the main theorems are given in Section 3. Finally, in Section 4, we discuss some possible extensions.
		
\section{Notation and main results}\label{results}
	In this section, we present our main results, including the strong consistency and the CLT of the urn composition. Before that, we introduce some needed notation.
	\subsection{Model and notations}
	We suppose $G=\{g_{0},g_{1}\dots, g_{d-1}\}$ is a finite Abelian group of order $d$, where $d$ is a positive integer. We denote by $\textgravedbl\ast\textacutedbl$ the binary multiplication on the group and write $z=a(x,y)=x\ast y$. Originally, we assume that
	the urn contains balls whose labels include arbitrary generators of group $G$. At stage $i$, we randomly draw $M_i$ balls separately with replacement, where $\{M_i,\ i=1,2,\dots\}$ is a sequence of independently and identically distributed random variables. If the labels of the $M_i$ balls are $X_{i1},\ X_{i2},\ \dots, \ X_{iM_i}$, then we add a ball to the urn with label $Z_i=X_{i1}\ast X_{i2}\ast\dots\ast X_{iM_i}$. After $n$ splits and generations, we define $\mathbf{S}_{n}=(S_{n,g_{0}},S_{n,g_{1}},\dots, S_{n,g_{d-1}})$ as the urn composition, where $S_{n,g_{i}}$ stands for the numbers of balls labeled $g_{i}$ after $n$ stages and $\sum_i S_{0,g_{i}}=n_{0}\geq1$. We note that $\{\mathbf{S}_{n}, n\ge1\}$ is a sequence of random $d$-vectors of nonnegative elements. We call this model a group multiplication add urn (GAU) model in the following.
		
	As a special case of the GAU model, if $G=({\mathbb Z}_{d},+)$, then $G=\{0,1,\dots,d-1\}$ and $z=a(x,y)=(x+y)\mod d$. With the rules of drawing and adding the balls in GAU model, the particular model is called a modular add urn (MAU) model. For brevity, we define $f(n)\preceq g(n)$ as $\lim\limits_{n\to \infty}f(n)/g(n)\le1$ and $f(n)\sim g(n)$ as $\lim\limits_{n\to \infty}f_{1}(n)/g_{1}(n)=1$. $\lfloor \cdot\rfloor$ stands for the floor function. Let $(a,b)$ be the greatest common divisor of two positive integers $a$ and $b$; let $b|a$ represent that $a$ is divisible by $b$. For all $g\in G$, $g^{-1}$ is defined as the inverse element of $g$ in the sequel.

	\subsection{Assumptions and main results} In this section, we present the first- and second-order convergences of the normalized urn composition of the Diaconis urn model with random multiple-drawing. We first deduce the first-order convergence of $\mathbf{S}_{n}$ under the MAU model, and generalize the result to the GAU model, then we derive the second-order convergence of the GAU model. Since the result is obvious when $d=1$, in the sequel, we assume that $d\ge2$. All the proofs of the results are postponed to Section \ref{proof}.

	Let ${\cal F}_n$ be the $\sigma$-field generated by $\{X_{i1},X_{i2},\dots, X_{iM_i},i=1,\dots, n\}$. The main results of this paper are obtained under the following assumptions.
	\begin{Assumption}\label{as1}
		At stage $i$, $i=1,2,\dots$, the $j$th ball $X_{ij}$ drawn is conditionally independent of the number of balls drawn $M_i$ given ${\cal F}_{i-1}$. Moreover, when given $M_i=m$, for all $j\in\{1,2,\dots, m\}$,
			$$X_{ij}|{\cal F}_{i-1}\sim \mbox{Multinomial}(1;\mathbf{S}_{i-1}/|\mathbf{S}_{i-1}|),$$
		where $|\mathbf{S}_{i-1}|=\sum\limits_{j=0}^{d-1}S_{i-1,g_j}$.	
	\end{Assumption}
	\begin{Assumption}\label{as2}
		Let $c\ge1$, $C\ge2$ be two finite integers. We define $m_1,m_2,\dots, m_c$ as all possible values of $M_i(i=1,2,\dots)$ that satisfy $2\le m_1<m_2<\dots<m_c\le C$. We assume that ${M_{i}\stackrel{i.i.d}\sim \mbox{Multinomial}(1;(p_1,p_2,\dots, p_c))}$. That is, for all $i=1,2,\dots$,
			$$P(M_i=m_1)=p_1,\ P(M_i=m_2)=p_2,\dots, \ P(M_i=m_c)=p_c.$$
		And, we assume that there exists at least one $s\in\{1,2,\dots, c\}$ such that $(m_s-1,d)=1$.
	\end{Assumption}
	\begin{remark}
		Here, we allow the number of balls drawn at each stage to be a random variable rather than a fixed number of two balls. The number of balls drawn, $M_i$, greater than or equal to 2 is the problem proposed by P. Diaconis.
	\end{remark}
		Assumption \ref{as2} ensures that all elements can be generated by a generator under the Diaconis urn model with random multiple-drawing. We illustrate the importance of this assumption based on a counterexample. We assume the group $G=(\mathbb{Z}_5,+)$, $c=2$, $m_1=6$, $m_2=11$, and the urn only contains balls initially labeled by the element $2$. Then, the urn still contains balls only labeled by $2$ before the second draw due to $d|(m_1-1)$ and $d|(m_2-1)$, which means that no matter how many times we draw and replace balls, the urn only contains balls labeled by $2$. Based on this fact, we give the following proposition.
	\begin{proposition}
			In the GAU model, if $G$ is a cyclic group, the urn initially contains balls only labeled by the generator of the group $G$. Under Assumption \ref{as2}, if $d|(m_s-1)$ for all $s\in\{1,2,\dots, c\}$, then no matter how many times we draw and replace balls, the urn still only contains balls labeled by this generator.
	\end{proposition}

	\begin{theorem}[Strong\ consistency\ of\ $\mathbf{S}_{n}$\ under\ MAU] \label{thm1}
		Under Assumptions \ref{as1}$\sim$\ref{as2}, for the MAU model, we have that for all $j\in\{0,1,\dots ,d-1\}$,
		$$a_{nj}\to 1/d,\ \ \ a.s.,$$
		where $a_{nj}=S_{nj}/(n+n_0)$.
	\end{theorem}
	\begin{remark}
		Theorem \ref{thm1} states the first-order convergence of the normalized numbers of balls in the urn, which implies that the sequence of the normalized urn composition $\mathbf{S}_{n}/n$ converges almost surely to a uniform distribution over group $G$ for the MAU model.
	\end{remark}
		Any Abelian group with prime order is a cyclic group and is isomorphic to (${\mathbb Z}_{d}$,+). Therefore, we have the following corollary.
		
	\begin{corollary} \label{coro1}
		Under Assumptions \ref{as1}$\sim$\ref{as2}, if $G$ is an Abelian group of order $d$, and $d$ is a prime number, then the sequence of the normalized urn composition converges to the uniform distribution over group $G$ almost surely.
	\end{corollary}	
		
	The following theorem generalizes Corollary \ref{coro1} to the case that $d$ is any positive integer.
	\begin{theorem}[Strong\ consistency\ of\ $\mathbf{S}_{n}$\ under\ GAU] \label{thm4}
		Under Assumptions \ref{as1}$\sim$\ref{as2}, for the GAU model, we have that for all $j\in\{0,1,\dots, d-1\}$,
		$$a_{n,g_{j}}\to 1/d, \ \ \ a.s.,$$
		\ where\ $a_{n,g_{j}}=S_{n,g_{j}}/(n+n_0)$.
	\end{theorem}
	\begin{remark}
		We consider the classical multicolor P\'{o}lya urn model with $\mathbf{S}_0=(1,1,\dots, 1)_{1\times d}$ and a ball which is randomly drawn from the urn and returned along with one ball of the same color at each stage. After $n$ stages, \cite{r34} showed that the proportion of balls of each color converges to a beta distribution $Beta(1,d-1)$ almost surely. We note that the mean of $Beta(1,d-1)$ is $1/d$, which equals the limiting value of the proportion under the GAU model. However, for the GAU model, even if there is a class of balls in the urn with a ratio very close to 1 at a certain moment, extracting multiple balls and performing group multiplications equalize the number of balls of each type in the urn when $n$ tends to infinity. When one ball is drawn with replacement at each stage, the ball with a higher initial percentage may be drawn with a higher probability at the next stage, resulting in an increasing number of balls in that category.
	\end{remark}
		
	Theorem \ref{thm4} is the result of the first-order convergence of $\mathbf{S}_{n}$. The next theorem states the second-order asymptotic result for $\mathbf{S}_{n}$.
	\begin{theorem}[Asymptotic\ normality\ of\ $\mathbf{S}_{n}$] \label{thm3}
		Under Assumptions \ref{as1}$\sim$\ref{as2}, for the GAU model, we have
		$$\frac{1}{\sqrt{n}}\bigg({{\bf S}_{n}}-\frac{n}{d}{\bf 1}_{d}\bigg)\to N_{d}(0,\Sigma)$$\ \ \ \it{in distribution,
			where}
		\begin{equation}  \label{115}
			\Sigma=
			\begin{pmatrix}
				\frac{d-1}{d^{2}} &  -\frac{1}{d^{2}} & \cdots& -\frac{1}{d^{2}} \\			
				-\frac{1}{d^{2}} & \frac{d-1}{d^{2}} & \cdots& -\frac{1}{d^{2}} \\
				\cdots & \cdots & \cdots & \cdots \\
				-\frac{1}{d^{2}} & -\frac{1}{d^{2}} & \cdots & \frac{d-1}{d^{2}} 
			\end{pmatrix}
			=\frac{1}{d}{\bf I}_d-\frac{1}{d^{2}}{\bf 1}_{d}{\bf 1}^{T}_{d},
		\end{equation}
		${\bf I}_d$ is the identity matrix of size $d$, ${\bf 1}_{d}$ is the $d$-dimensional all-ones vector, and ${\bf 1}^{T}_{d}$ denotes the transpose of ${\bf 1}_{d}$.
	\end{theorem}
	\begin{remark}				
		We suppose that $\mathbf{W}_{i}$, $i=1,2,\dots$, are independent and identically distributed $\rm Multinomial\ (1;1/\it d, \dots,  1/\it d)$ random variables. Then, $(\sum_{i=1}^{n}\mathbf{W}_{i}-n{\bf 1}_{d}/d)/\sqrt{n}$ weakly converges to $N_{d}(0,\Sigma)$, where $\Sigma$ is defined in \eqref{115}.
		In fact, for the GAU model, we can imagine from Theorem \ref{thm4} that when $n$ is large enough, the proportion of each type of ball in the urn is close to $1/d$. Therefore, when ${\mathcal{F}}_{n}$ is given, the result of randomly drawing a ball at the $(n+1) th$ stage has approximately the same distribution as $\mathbf{W}_{n}$. This fact leads to the distribution of the ball put in the urn at the $(n+1)th$ stage also approximately obeying the same distribution as $\mathbf{W}_{n}$ due to the special addition rules. Thus, the process of drawing and replacing the ball can be approximately regarded as independent; as a result, $(\sum_{i=1}^{n}\mathbf{W}_{i}-n{\bf 1}_{d}/d)/\sqrt{n}$ and $(\mathbf{S}_{n}-n{\bf 1}_{d}/d)/\sqrt{n}$ have the same asymptotic distribution.
	\end{remark}

	\begin{remark}
		If $c=1$, $m_1=2$ and $P(M_i=m_1)=P(M_i=2)=1$ for all $i\in\{1,2,\dots\}$, then $(m_1-1,d)=1$. Thus, for Abelian groups, our results indicate that the first- and second-order convergence of the normalized urn composition hold for the model proposed by Diaconis.
	\end{remark}
		
\section{Technical proofs}\label{proof}
	In this section, we present the technical proofs of the main results in Section \ref{results}.
	In Assumption \ref{as2}, we assume that there exists at least one $i\in\{1,2,\dots, c\}$ such that $(m_i-1,d)=1$. Without loss of generality, we assume that $(m_1-1,d)=1$.

	\subsection{Proof of Theorem \ref{thm1}}\label{sec3.1}
	The strategy of this proof is as follows. Since the urn contains some balls whose labels include the generators of group $G$, the main task of this proof is to show that the number of balls labeled with each element of $G$ in the urn is of the same order as $n$ almost surely, i.e.,
	\begin{equation}\label{go1}
		\limsup_{n\to\infty}\min_{j}a_{nj}>0,\ \ a.s..
	\end{equation}
	When $d$ is a prime number, for all $j\in\{1,2,\dots, d-1\}$, $j$ is a generator of group $G$. Therefore, to prove (\ref{go1}), it is critical to show that there exists $ j\in\{1,2,\dots, d-1\}$, such that
	\begin{equation}\label{go2}
		\limsup_{n\to\infty}a_{nj}>0,\ \ a.s.,
	\end{equation}
	which is equivalent to
	\begin{equation}\label{go3}
		\liminf_{n\to\infty}a_{n0}<1,\ \ a.s..
	\end{equation}			
	More specifically, we divide the proof procedure into the following steps. Defining $\widetilde{C}=\max\left\{\sqrt[m_1-1]{\frac{3}{2m_1}},\sqrt[m_2-1]{\frac{3}{2m_2}},\dots,  \sqrt[m_c-1]{\frac{3}{2m_c}}\right\}$ and ${\mathcal{G}}=\liminf\limits_{n\to\infty}\left\{a_{n0}>\widetilde{C}\right\}$, we have
	\begin{enumerate}[label=Step \Roman*:]
		\item $\lim\limits_{n\to\infty} \Big\{\sum\limits_{j\ge 1}S_{nj}>k\Big\}$,\ \ a.s.,\ $\forall$ $k\geq0$.
		\item On the set
		${\cal G}$,\ except for a null set, we have
		\begin{equation}
			\liminf\limits_{n\to\infty} \Big\{\sum\limits_{j\ge1}S_{nj}/(\ln
			n)^k\ge N\Big\},\ \forall\ k\geq0,\ N>0.  \nonumber
		\end{equation}
		\item On the set
		${\cal G}$,\ except for a null set, we have
		\begin{equation}
			\liminf\limits_{n\to\infty}
			\Big\{\sum\limits_{j\ge1}S_{nj}/n^{\tau}\ge N\Big\},\ \exists \tau>3/4, \ \forall \ N>0.  \nonumber
		\end{equation}
		\item $\liminf\limits_{n\to\infty} a_{n0}<1$,\ \ a.s.
		\item If $d$ is a prime number, then $a_{nj}\to 1/d$,\ \ a.s.
		\item Generalize the result to any $d$.
	\end{enumerate}
	It is worth noting that the conclusions in Steps II and III are obvious on the set ${\mathcal {G}}^{c}$ from the fact that  $\mathcal {G}\supset A:=\left\{\liminf\limits_{n\to\infty} a_{n0}>\widetilde{C}\right\}$. In the following, we prove Theorem \ref{thm1} step by step.

	\subsubsection{Proof of Step I}

	Before proving Step I, we state the following results.
	\begin{lemma} \label{le1}
	We assume that the urn contains only the balls labeled as the generator $q$ of the group $G$, and Assumptions \ref{as1}$\sim$\ref{as2} hold. We define the following two stopping times:
	\begin{eqnarray*}
	&&T_1:=\begin{cases}\inf \{n\ge1:M_n=m_1\},&\text{if\ \ $ \{n\ge1:M_n=m_1\}\neq\varnothing;$}\\
	\infty,\ &\text{otherwise,} \end{cases}\\
	&&T_2:=\begin{cases} \inf\{n\ge1:&\mbox{exists}\ j\in\{0,1,\dots, q-1,q+1,\dots, d-1\}\ \mbox{such\ that\ } a_{n,j}>0\},\\
		 &\text{if\ \ $\{n\ge1:\exists\ j\in\{0,1,\dots, q-1,q+1,\dots, d-1\}\ \mbox{s.t.\ }a_{n,j}>0\}\neq\varnothing;$}\\
	\infty,\ &\text{otherwise.} \end{cases}
	\end{eqnarray*}
	Then, we have
	$$T_2<\infty\ \ \ a.s..$$
	\end{lemma}
	The proof of Lemma \ref{le1} is postponed to Section \ref{prfles}. Lemma \ref{le1} indicates the following facts. We assume that the urn initially contains balls labeled by only one kind of element of group $G$, and this element is the generator of the group. Then, under Assumptions \ref{as1}$\sim$\ref{as2}, another kind of element must be generated after sufficiently many finite stages.

	In the following, we prove Step I by mathematical induction. Since the urn initially contains balls whose labels include the generators of group $G$, $\sum_{j\geq1}S_{1j}\geq1$, which implies Step I holds for $k=0$. We assume that $k\geq1$ and Step I holds for $k-1$. We define the set
		$$
{\cal A}_n(k)=\bigg\{\sum_{j\ge 1}S_{nj}\le k\bigg\}.
		$$
	It is obvious that ${\cal A}_n(k)$ monotonically decreases as $n\to\infty$. Thus, ${\cal A}(k)=\lim\limits_{n\to\infty}{\cal A}_n(k)$ exists.
	Next, we prove that
	\begin{equation} \label{gong12}
		P({\cal A}(k))=0.
	\end{equation}
	It follows that
	\begin{eqnarray}
		&&P({\cal A}_{n+1}(k))\nonumber\\
		&=&
		P({\cal A}_n(k-1))+P\left(\left\{\sum_{j\ge 1}S_{n,j}=k\right\}\cap\left\{Z_{n+1}=0\right\}\right)\nonumber\\
		&=&P({\cal A}_n(k-1))+\sum\limits_{s=1}^{c}\mathbb{E}\left[I\left(\sum_{j\ge 1}S_{n,j}= k\right)P(Z_{n+1}=0,M_{n+1}=m_s|{\cal F}_n)\right]\nonumber\\
		&=&P({\cal A}_n(k-1))+\sum\limits_{s=1}^{c}\mathbb{E}\left[I\left(\sum_{j\ge 1}S_{n,j}= k\right)P(Z_{n+1}=0|M_{n+1}=m_s,{\cal F}_n)p_s\right]\nonumber\\
		&=&P({\cal A}_n(k-1))+\sum\limits_{s=1}^{c}p_s\mathbb{E}\left[I\left(\sum_{j\ge 1}S_{n,j}= k\right)\sum\limits_{x_{1}\ast x_{2}\ast\dots\ast x_{m_s}=0}a_{n,x_1}a_{n,x_2}\dots a_{n,x_{m_s}}\right]\nonumber\\
		&\le&P({\cal A}_n(k-1))+\mathbb{E}\left[I\left(\sum_{j\ge 1}S_{n,j}= k\right)\left(1-\frac1{n+n_0}\right)\right]\nonumber\\
		&\le&\frac1{n+n_0}P({\cal A}_n(k-1))+P({\cal A}_n(k))\left(1-\frac1{n+n_0}\right).\label{gong22}
	\end{eqnarray}
	Here, we used the fact that $a_{n,d}=a_{n,0}$ and that for all $s\in\{1,2,\dots, c\}$,
	\begin{eqnarray*}
		&&1-\sum\limits_{x_{1}\ast x_{2}\ast\dots\ast x_{m_s}=0}a_{n,x_1}a_{n,x_2}\dots a_{n,x_{m_s}}\\
		&=&\sum\limits_{x_{1}\ast x_{2}\ast\dots\ast x_{m_s}=0}a_{n,x_1}a_{n,x_2}\dots a_{n,x_{m_s-1}}-\sum\limits_{x_{1}\ast x_{2}\ast\dots\ast x_{m_s}=0}a_{n,x_1}a_{n,x_2}\dots a_{n,x_{m_s}}\\
		&=&\sum\limits_{x_{1}\ast x_{2}\ast\dots\ast x_{m_s}=0}a_{n,x_1}a_{n,x_2}\dots a_{n,x_{m_s-1}}(1-a_{n,x_{m_s}})\\
		&\ge&\min\limits_{j}(1-a_{n,x_j})\ge\frac{1}{n+n_0},
	\end{eqnarray*}
	where the last inequality is because the urn contains at least two kinds of balls after sufficiently many stages by Lemma \ref{le1}. In the first equality, we use that for all $j\in\{0,1,\dots, d-1\}$,
	\begin{eqnarray}\label{go949}
		&&\sum\limits_{x_1\ast x_2\ast\dots\ast x_{m_s}=j}a_{n_k,x_1}a_{n_k,x_2}\dots a_{n_k,x_{m_s-1}}\\
		&=&\sum_{x_{m_s}=0}^{d-1}\sum\limits_{x_1\ast x_2\ast\dots\ast x_{m_s-1}=x^{-1}_{m_s}\ast j}a_{n_k,x_1}a_{n_k,x_2}\dots a_{n_k,x_{m_s-1}}\nonumber\\
		&=&\sum\limits_{x_1\ast x_2\ast\dots\ast x_{m_s-1}=0}^{d-1}a_{n_k,x_1}a_{n_k,x_2}\dots a_{n_k,x_{m_s-1}}=1,\nonumber
	\end{eqnarray}

	Recursively employing inequality (\ref{gong22}), we obtain that
	\begin{eqnarray*}
		P({\cal A}_{2n}(k))&\le& \sum_{i=n}^{2n-1}\frac1{i+n_0}P({\cal A}_{i}(k-1))+\frac{n+n_0-1}{2n+n_0-1} P({\cal A}_{n}(k))\\
		&\le&P({\cal A}_n(k-1))+\frac{n+n_0-1}{2n+n_0-1}P({\cal A}_n(k)).
	\end{eqnarray*}
	Let $n\to \infty$; by induction, we obtain $P({\cal A}_n(k-1))\to 0$, which implies ({\ref{gong12}}). Then, the proof of Step I is complete.
		
\subsubsection{Proof of Step II}
	Let $p_n=S_{n1}+\cdots+S_{n,d-1}$ and $q_{j,n}=p_n-p_{j}$.

	We recall that ${\cal F}_n$ is the $\sigma$-field generated by $\{(X_{i1},X_{i2},\dots, X_{iM_i}),i=1,\dots, n\}$, where $M_i$ denotes the number of balls drawn at stage $i$. Then, we have that
	\begin{eqnarray*}
		P(Z_{i}\neq 0|{\cal F}_{i-1})&=&\sum\limits_{s=1}^{c}P(Z_{i}\neq 0,M_i=m_s|{\cal F}_{i-1})=\sum\limits_{s=1}^{c}p_sP(Z_{i}\neq 0|M_i=m_s,{\cal F}_{i-1})\\
		&\le&\sum\limits_{s=1}^{c}p_s\left[P(X_{i1}\neq 0|{\cal F}_{i-1})+P(X_{i2}\neq 0|{\cal F}_{i-1})+\dots+P(X_{im_s}\neq 0|{\cal F}_{i-1})\right]\\
		&=&\sum\limits_{s=1}^{c}p_sm_s(a_{i-1,1}+a_{i-1,2}+\dots+a_{i-1,d-1})\le\frac{Cp_{i-1}}{i-1+n_0}\le\frac{Cp_{i-1}}{i}.
	\end{eqnarray*}
	On the other hand, if $a_{i-1,0}>\widetilde{C}$ and $i\ge 2(n_0-1)$, we have
	\begin{eqnarray*}
		P(Z_{i}=0|{\cal F}_{i-1})&=&\sum\limits_{s=1}^{c}p_sP(Z_{i}=0|M_i=m_s,{\cal F}_{i-1})\\
		&\le&\sum\limits_{s=1}^{c}p_s[1-P(\{X_{i1}\neq 0\}\cap\{X_{i2}=0\}\cap\{X_{i3}=0\}\cap\dots\cap\{X_{im_s}=0\}|{\cal F}_{i-1})\\
		&&-P(\{X_{i1}=0\}\cap\{X_{i2}\neq 0\}\cap\{X_{i3}=0\}\cap\dots\cap\{X_{im_s}=0\}|{\cal F}_{i-1})-\dots\\
		&&-P(\{X_{i1}=0\}\cap\{X_{i2}=0\}\cap\{X_{i3}\neq 0\}\cap\dots\cap\{X_{im_s}\neq0\}|{\cal F}_{i-1})]\\
		&=&\sum\limits_{s=1}^{c}p_s\left[1-m_sa^{m_s-1}_{i-1,0}(a_{i-1,1}+\dots+a_{i-1,d-1})\right]<1-\frac{p_{i-1}}{i}.
	\end{eqnarray*}
	Let
		$$
	{\cal G}_{j,n}=\bigcap_{\ell=j}^n\left\{a_{\ell,0}>\widetilde{C}\right\}.
		$$
	Then, for $\delta<1$ and $\lfloor n^{\delta}\rfloor\geq2(n_{0}-1)$,
	\begin{eqnarray*}
		&&P\left({\cal G}_{\lfloor n^{\delta}\rfloor,n-1}\cap \{ q_{\lfloor n^{\delta}\rfloor,n}=k\}|{\cal F}_{\lfloor n^{\delta}\rfloor}\right)\\
		&=& \sum_{\lfloor n^{\delta}\rfloor<j_1<\cdots<j_k\le n}P\left.\left({\cal G}_{\lfloor n^{\delta}\rfloor,n-1}\cap \left\{{Z_j\neq 0,j= j_1,\cdots,j_k;\atop Z_j=0,j\neq j_1,\cdots,j_k }\right\}\right|{\cal F}_{\lfloor n^{\delta}\rfloor}\right)\\
		&\le &\sum_{\lfloor n^{\delta}\rfloor<j_1<\cdots<j_k\le n}C^k\prod_{{\lfloor n^{\delta}\rfloor}<j<j_1}
		\left(1-\frac{p_{\lfloor n^{\delta}\rfloor}}{j}\right)\prod_{t=1}^{k}\left[\frac{p_{\lfloor n^{\delta}\rfloor}+t-1}{j_t}\prod_{j_t<j<j_{t+1}}
		\left(1-\frac{p_{\lfloor n^{\delta}\rfloor}+t}{j}\right)\right],
	\end{eqnarray*}
	where $j_{k+1}=n+1$. By the identity
		$$
	\frac{a}{b}\left(1-\frac{a+1}{b+1}\right)=\left(1-\frac{a}{b}\right)\frac{a}{b+1},
		$$
	we have that
	\begin{eqnarray*}
		&&P\left({\cal G}_{\lfloor n^{\delta}\rfloor,n-1}\cap \{q_{\lfloor n^{\delta}\rfloor,n}=k\}|{\cal F}_{\lfloor n^{\delta}\rfloor}\right)\\
		&\le &C^k\sum_{\lfloor n^{\delta}\rfloor<j_1<\cdots<j_k\le n}\prod_{\lfloor n^{\delta}\rfloor<j\le n-k}\left(1-\frac{p_{\lfloor n^{\delta}\rfloor}}{j}\right)\prod_{t=1}^k\frac{p_{\lfloor n^{\delta}\rfloor}+t-1}{n-k+t}\\
		&=&C^k{n-\lfloor n^{\delta}\rfloor \choose k}\prod_{\lfloor n^{\delta}\rfloor<j\le n-k}\left(1-\frac{p_{\lfloor n^{\delta}\rfloor}}{j}\right)\prod_{t=1}^k
		\frac{p_{\lfloor n^{\delta}\rfloor}+t-1}{n-k+t}\\
		&\le&C^k\prod_{\lfloor n^{\delta}\rfloor<j\le n-k}\left(1-\frac{p_{\lfloor n^{\delta}\rfloor}}{j}\right)\prod_{t=1}^k
		\frac{p_{\lfloor n^{\delta}\rfloor}+t-1}{t},
	\end{eqnarray*}
	which implies
	\begin{eqnarray} \label{key1}
		&&P\left({\cal G}_{\lfloor n^{\delta}\rfloor,n-1}\cap \{q_{\lfloor n^{\delta}\rfloor,n}\le K\}|{\cal F}_{\lfloor n^{\delta}\rfloor}\right)\nonumber\\
		& \le & \sum_{k=0}^KC^k\prod_{\lfloor n^{\delta}\rfloor<j\le n-k} \left(1-\frac{p_{\lfloor n^{\delta}\rfloor}}{j}\right) \prod_{t=1}^k\frac{p_{\lfloor n^{\delta}\rfloor}+t-1}{t}\nonumber\\
		& \le & C^{K}\prod_{\lfloor n^{\delta}\rfloor<j\le n-K} \left(1-\frac{p_{\lfloor n^{\delta}\rfloor}}{j}\right) {p_{\lfloor n^{\delta}\rfloor}+K \choose K}\nonumber\\
		& \le & C^{K}\exp\left[\ln \left(1+\frac{p_{\lfloor n^{\delta}\rfloor}}{K}\right)^{K}+\ln\left(1+\frac{K}{p_{\lfloor n^{\delta}\rfloor}}\right)^{p_{\lfloor n^{\delta}\rfloor}}- \ln \left(\frac{n-K}{\lfloor n^{\delta}\rfloor}\right)^{p_{\lfloor n^{\delta}\rfloor}}\right]\nonumber\nonumber\\
		& \preceq & \exp\{-p_{\lfloor n^{\delta}\rfloor}[(1-\delta)\ln n-2]+K(1+\ln C)\}.
	\end{eqnarray}
	In the first step, we use the fact that if $q+1\le n$, then
	$${n \choose q}+{n\choose q+1}={n+1 \choose {q+1}},$$
	and in the last step, we assume that $K=o(n)$.
		
	To prove step II, it is enough to show that
	\begin{equation}
		P\left({\cal G}\cap \{p_n< N(\ln n)^k,\   i.o.\}\right)=0, \ \forall k\ge1,\ N\ge \lfloor C^{2}-C-\ln C-1 \rfloor+1.   \label{key2}
	\end{equation}
	We show (\ref{key2}) by induction.
	We recall that ${\mathcal{G}}=\liminf\limits_{n\to\infty}\left\{a_{n0}>\widetilde{C}\right\}$. From Step I, we know that (\ref{key2})
	holds for $k=0$ and any $N>1$. For any $k\ge 1$ if
		$$
	P\left({\cal G}\cap \{p_n<NC^{k+1}(\ln n)^{k-1},\ i.o.\}\right)=0.
		$$
	Then, it follows from (\ref{key1}) that
	\begin{eqnarray*}
		&&P\left({\cal G}\cap \{p_n< N(\ln n)^k,\  i.o.\}\right)\\
		&\le& P\left({\cal G}\cap \{p_{n}< NC^{k+1}(\ln n)^{k-1},\  i.o.\}\right)\\
		&&+P\left({\cal G}\cap \left\{\{p_n< N(\ln n)^k\}\cap \{p_{\lfloor n^{1/C}\rfloor}\ge NC^{k+1}(\ln n^{1/m})^{k-1}\},\ i.o.\right\}\right)\\
		&\le& P\left({\cal G}\cap \{p_n< NC^{k+1}(\ln n)^{k-1},\  i.o.\}\right)\\
		&&+\lim_{t\to\infty}\sum_{n=t}^{\infty} P\left({\cal G}_{\lfloor n^{1/C}\rfloor,n-1}\cap\left\{ \{p_n< N(\ln n)^k\} \cap \{p_{\lfloor n^{1/C}\rfloor}\ge C^{2}N(\ln n)^{k-1}\}\right\}\right)\\
		&=&0. 
	\end{eqnarray*}
	Here, we use the fact that
	\begin{eqnarray*}
	&&P\left({\cal G}_{\lfloor n^{1/C}\rfloor,n-1}\cap\left\{ \{p_n< N(\ln n)^k\}  \cap \{p_{\lfloor n^{1/C}\rfloor}\ge C^{2}N(\ln  n)^{k-1}\}\right\}\right)\\
	&\preceq&\exp\{-N[C^{2}-C-\ln C-1](\ln n)^{k}\},
	\end{eqnarray*}
	which is summable since $N\ge \lfloor1/(C^2-C-\ln C-1)\rfloor+1$. This completes the proof of Step II.

\subsubsection{Proof of step III}
	Now, we again use (\ref{key2}) to further improve the convergence rate of $p_n$.
	For a fixed, large $n$, we choose $\delta=\delta_n=1-A/\ln n$, where $A>0$ is a large constant determined later. Let
	\begin{eqnarray*}
		n_j&=&\lfloor n^{\delta^{j}}\rfloor\\
		K_j&=&\frac{B^{l_n-j}\delta^{-(l_n-j)(l_n-j+1)/2}}{2^{l_n-j}(1+\ln C)^{l_n-j+1}}(B/A)^{\alpha}B(\ln n)^{\alpha},
	\end{eqnarray*}
	where $j=0,1,\dots, l_n=\lfloor A^{-1}\ln n(\ln A-\ln B)\rfloor$ with $B<A$ and $\alpha$ is a constant with $2^{\alpha-1}\delta^{\alpha}\ge 2$. We note that $K_j=\frac{B\delta^{j-l_n}}{2(1+\ln C)}K_{j+1}$, $\ln n_{l_n+1}=\delta\ln n_{l_n}$ and
	\begin{eqnarray}
		(1-\delta)\ln n_{l_n}\sim (1-\delta)\delta^{l_n}\ln n=A e^{l_n\ln \delta}\sim B . \label{kkey0}
	\end{eqnarray}
	Thus, we have that $\ln n_{l_n} \sim BA^{-1}\ln n$.
	Then, by Step II, we may assume,
	$p_{n_{l_n+1}}> 2^\alpha(\ln n_{l_n+1})^\alpha$, and hence
	\begin{eqnarray*}
		&&p_{n_{l_n+1}}[(1-\delta)\ln n_{l_n}-2]-K_{l_n}(1+\ln C)\\
		&\ge & (2B/A)^{\alpha} B(\ln n)^{\alpha}\delta^\alpha/2-K_{l_n}(1+\ln C),\\
		&\ge & 2(B/A)^{\alpha}B(\ln n)^{\alpha}-(B/A)^{\alpha}B(\ln n)^{\alpha}\\
		&=&(B/A)^{\alpha}B(\ln n)^{\alpha}. 
	\end{eqnarray*}
	Applying (\ref{key1}) with $n=n_{l_n}$, we have
	\begin{eqnarray}
		&&P\left({\cal G}_{n_{l_n+1},n_{l_n}-1}\cap \left\{q_{n_{l_n+1},n_{l_n}}
		\le K_{l_n}\right\}|{\cal F}_{n_{l_n+1}}\right)\nonumber\\
		&\preceq& \exp\left(-(B/A)^{\alpha}B(\ln n)^{\alpha}\right)  \label{kkey1}
	\end{eqnarray}
	For $j<l_n$, we have
		$$
(1-\delta)\ln n_j\sim \delta^{j}(1-\delta)\ln n\sim B\delta^{j-l_n}.
		$$
If ${q_{n_{j+2},n_{j+1}}>K_{j+1}}$, then $p_{n_{j+1}}> K_{j+1}$. Again, by (\ref{key1}), we have that
	\begin{eqnarray}
		&&P\left({\cal G}_{n_{j+1},n_j-1}\cap\left\{q_{n_{j+2},n_{j+1}}> K_{j+1}, q_{n_{j+1},n_{j}}\le K_{j}\right\}|{\cal F}_{n_{j+1}}\right)\nonumber\\
		&\preceq&\exp\left(-K_{j+1}[(1-\delta)\ln n_j-2]+K_j(1+\ln C) \right)\nonumber\\
		&\le&\exp\left(-(B\delta^{j-l_n}-2)K_{j+1}+K_j(1+\ln C)\right)\nonumber\\
		&\le&\exp\left(-K_{j}\right)\nonumber\\
		&=&\exp\left(-\frac{B^{l_n-j}\delta^{-(l_n-j)(l_n-j+1)/2}}{2^{l_n-j}(1+\ln C)^{l_n-j+1}}B(BA^{-1}\ln n)^{\alpha}\right).   \label{kkey2}
	\end{eqnarray}
	By (\ref{kkey1}) and (\ref{kkey2}), we find
	\begin{eqnarray*}
		&&P\left({\cal G}_{n_{l_n+1},n-1}\bigcap\left.\left\{\bigcup_{j=0}^{l_n}\left\{q_{n_{j+1},n_{j}}\le K_{j}\right\}\right\}\right|{\cal F}_{n_{l_n+1}}\right)\\
		&=&P\left({\cal G}_{n_{l_n+1},n-1}\bigcap\left.\left\{\bigcup_{j=0}^{l_n-1}\left\{q_{n_{j+2},n_{j+1}}>K_{j+1},q_{n_{j+1},n_{j}}\le K_{j}\right\}\bigcup\left\{q_{n_{l_{n}+1},n_{l_{n}}}\le K_{l_{n}}\right\}\right\}\right|{\cal F}_{n_{l_n+1}}\right)\\
		&\le& \sum_{j=0}^{l_n}\exp\left(-\frac{B^{l_n-j}\delta^{-(l_n-j)(l_n-j+1)/2}}{2^{l_n-j}(1+\ln C)^{l_n-j+1}}B(BA^{-1}\ln n)^{\alpha}\right),
	\end{eqnarray*}
	which is summable. By the Borel-Cantelli lemma, when ${\cal G}$ occurs, except on the null set, for all large $n$, we have that
		$$
	p_n > K_0.
		$$

	Note that
		$$
	K_0\sim \frac{B^{l_n}\delta^{-l_n(l_n+1)/2}}{2^{l_n}(1+\ln C)^{l_n+1}}B(BA^{-1}\ln n)^\alpha.
		$$
	From (\ref{kkey0}), we have that
		$$
	\delta^{-l_n}=AB^{-1}.
		$$
	Therefore, for all large $n$, $K_0> \eta^{\ln n}(BA^{-1}\ln
	n)^\alpha$, where $\eta=\left(\sqrt{AB}/(2(1+\ln C))\right)^{A^{-1}\ln(A/B)}$. Choosing $A>B$ large so that $\eta>1$, we have proved that when ${\cal G}$ occurs, except for a null set, for
	all large $n$,
	\begin{equation} p_n>n^{\tau_0}(BA^{-1}\ln n)^\alpha,
	\end{equation}
	where $\tau_0=\ln \eta$.
		
	If $\tau_0<3/4$, then by repeating the above arguments with
	$p_{n_{l_n+1}}\ge M n_{l_n}^{\tau_0}=M n^{\rho\tau_0}$ for any
	$M>1$, we can prove that when ${\cal G}$ occurs, except for as a null set, for
	all large $n$,
	$$p_n>M n^{\tau_0}n^{\rho\tau_0}.$$
	Note that
		$$ 
	n^{\tau_0}n^{\rho\tau_0}=n^{(1+\rho)\tau_0},
		$$
	where $\rho=1-\tau_0=B/A$. Here, we can choose $A$ and $B$ large and
	$B/A=\rho$. Continuing this procedure, we have that for any
	$k>0$,
		$$
	p_n>Mn^{\tau_0}n^{\rho\tau_0}\cdots n^{\rho^{k}\tau_0}=Mn^{(1+\rho+\cdots+\rho^k)\tau_0}.
		$$
	Then, by the fact that $k\to\infty$, $(1+\rho+\cdots+\rho^k)\tau_0\to 1$, we complete the proof of Step III.
		
\subsubsection{Proof of Step IV}
	We recall that $A=\{\liminf\limits_{n\to\infty}a_{n0}>\widetilde{C}\}$; we need to prove Step IV only on the set $A$ because Step IV holds obviously on $A^{c}$.	
	Let $n_1=1$ and $n_k=n_{k-1}+k=\frac12k(k+1)$ for $k\ge 2$, and
		$$
	Q_{ni}=I(Z_n=i)-P(Z_n=i|{\cal F}_{n-1}),
		$$
	for $i=0,1,\cdots,d-1$. It follows from Bernstein's inequality (see Proposition 2.1 of \cite{r29}) that
	\begin{equation}\label{Bernine}
		P(\max_{n_k<n\le n_{k+1}}|\sum_{t=n_k+1}^{n} Q_{ti}|> \sqrt{k}\log k)\leq 2\exp\left\{-\frac{k\log^{2}k}{2\sqrt{k}\log k+k+1}\right\}.
	\end{equation}
	Here, we choose $\tau=n_{k+1}$, $b=k+1$ and $a=\sqrt{k}\log k$ in Proposition 2.1 of \cite{r29}.				
		
	Denote
		$$
	\Omega_0=\liminf\limits_{k\to\infty}\left\{\max_{n_{k-1}<n\le n_k}\left|\frac1{n+n_0}\sum_{t=n_{k-1}+1}^n Q_{tj}\right| \le 2k^{-3/2}\ln k\right\}.
		$$
	By \eqref{Bernine}, we know that
	\begin{equation}\label{gong34}
		P(\Omega_0)=1.
	\end{equation}
		
	Then, we show that on the set $A$, the following event holds:
	\begin{equation}\label{eqq97}
		\liminf\limits_{n\to\infty}\Big\{\sum\limits_{j\neq i}S_{nj}/n^{\tau}\ge N\Big\},\ \   \exists \tau>3/4,\  \forall i\in\{0,1,\dots, d-1\},\ \forall N>0.
	\end{equation}
	When $i=0$, we can verify that (\ref{eqq97}) holds from Step III directly. When $i\neq0$, (\ref{eqq97}) holds obviously because $A=\{\liminf\limits_{n\to\infty}a_{n0}>\widetilde{C}\}$. Then, we have that for all $N>0$, $m_s,\ s=1,2,\dots, c$ and a sufficiently large $n$,
	\begin{eqnarray} 
		&&\sum\limits_{x_{1}\ast x_{2}\ast\dots\ast x_{m_s}=0,\ x_{m_s}\neq 0}a_{n,x_1}a_{n,x_2}\dots a_{n,x_{m_s-1}} \label{gong29}\\
		&=&\sum_{x_2=0}^{d-1}\dots\sum_{x_{m_s-1}=0}^{d-1}\left(\sum_{x_{m_s}=1}^{d-1}\sum_{x_1\ast x_2\ast\dots x_{m_s}=0}a_{n,x_1}\right)a_{n,x_2}a_{n,x_3}\dots a_{n,x_{m_s-1}}\nonumber\\
		&\ge&\frac{N}{2}n^{\tau-1}\sum_{x_2=0}^{d-1}\dots\sum_{x_{m_s-1}=0}^{d-1}a_{n,x_2}a_{n,x_3}\dots a_{n,x_{m_s-1}}=\frac{N}{2}n^{\tau-1},\nonumber
	\end{eqnarray}
		
	For $n_{k-1}\le n\le n_k$, we have that for all large $k$,
	\begin{equation} \label{goxin}
		|a_{nj}-a_{n_{k-1},j}|=\left|\frac{S_{nj}-S_{n_{k-1},j}}{n+n_{0}}-\frac{(n-n_{k-1})S_{n_{k-1},j}}{(n+n_{0})(n_{k-1}+n_{0})}\right|<4/k,
	\end{equation}
	which implies that		
	\begin{align}
		&\frac1{n+n_0}\left|\sum_{i=n_{k-1}+1}^n[P(Z_i=0|{\cal F}_{i-1})-P(Z_{n_{k-1}+1}=0|{\cal F}_{n_{k-1}})]\right| \label{oo2} \\
		=&\frac1{n+n_0}\left|\sum_{i=n_{k-1}+1}^n\sum\limits_{s=1}^{c}p_s[P(Z_i=0|M_i=m_s,{\cal F}_{i-1})-P(Z_{n_{k-1}+1}=0|M_{n_{k-1}+1}=m_s,{\cal F}_{n_{k-1}})]\right|\nonumber \\
		=&\frac1{n+n_0}\left|\sum_{i=n_{k-1}+1}^n\sum\limits_{s=1}^{c}p_s\sum_{x_{1}\ast 	x_{2}\ast\dots\ast x_{m_s}=0}\left(a_{i-1,x_{1}}  \dots a_{i-1,x_{m_s}}-a_{n_{k-1},x_1}\dots a_{n_{k-1},x_{m_s}}\right)\right| \nonumber \\
		=&\frac1{n+n_0}\left|\sum_{i=n_{k-1}+1}^n\sum\limits_{s=1}^{c}p_s\sum_{x_{1}\ast 	x_{2}\ast\dots\ast x_{m_s}=0}\left(a_{i-1,x_{1}}  \dots a_{i-1,x_{m_s}}-a_{n_{k-1},x_1} \dots a_{n_{k-1},x_{m_s}}\right)\right| \nonumber \\
		\leq&\frac{8Cd^{C-1}}{k^2}. \nonumber
	\end{align}
	Recalling the definition of $n_k$, we have
	\begin{eqnarray*}
		a_{n_k,0}&=&\frac{S_{0,0}+\sum\limits_{i=1}^{n_k}I(Z_i=0)}{n_k+n_{0}}=\frac{n_{k-1}+n_0}{n_k+n_{0}}a_{n_{k-1},0}+\frac{1}{n_{k}+n_{0}}\sum_{i=n_{k-1}+1}^{n_k}I(Z_i=0).
	\end{eqnarray*}
	Therefore, on the set $\Omega_{0}\cap A$, if $a_{n_{k-1},0}>2/3$,
	\begin{align}\label{oo1}
		&a_{n_k,0}-a_{n_{k-1},0} \\
		=&\frac{1}{n_{k}+n_0}\sum_{i=n_{k-1}+1}^{n_k}[I(Z_i=0)-a_{n_{k-1},0}]  \nonumber\\
		=&\frac{1}{n_{k}+n_0}\sum_{i=n_{k-1}+1}^{n_k}Q_{i,0} +\frac{1}{n_{k}+n_0}\sum_{i=n_{k-1}+1}^{n_k}[P(Z_i=0|{\cal F}_{i-1})-P(Z_{n_{k-1}+1}=0|{\cal F}_{n_{k-1}})]  \nonumber \\
		&+\frac{1}{n_{k}+n_0}\sum_{i=n_{k-1}+1}^{n_k}[P(Z_{n_{k-1}+1}=0|{\cal F}_{n_{k-1}})-a_{n_{k-1},0}]  \nonumber \\
		\le&\frac{2\ln k}{k^{3/2}}+\frac{8Cd^{C-1}}{k^2}+\frac{k}{n_k+n_0}\left(\sum\limits_{s=1}^{c}p_s\sum_{x_{1}\ast x_{2}\ast\dots\ast x_{m_s}=0}a_{n_{k-1},x_1} \dots a_{n_{k-1},x_{m_s}}-a_{n_{k-1},0}\right)\nonumber \\
		\le&\frac{2\ln k}{k^{3/2}}+\frac{8Cd^{C-1}}{k^2}\nonumber\\
		&+\frac{k}{n_k+n_0}\sum\limits_{s=1}^{c}p_s\sum_{x_{1}\ast x_{2}\ast\dots\ast x_{m_s}=0}a_{n_{k-1},x_1}\dots a_{n_{k-1},x_{m_s-1}}(a_{n_{k-1},x_{m_s}}-a_{n_{k-1},0})\nonumber \\
		\le&\frac{3\ln k}{k^{3/2}}-\frac{k}{3n_k+3n_0}\sum\limits_{s=1}^{c}p_s
		\sum_{x_{1}\ast x_{2}\ast\dots\ast x_{m_s}=0, x_{m_s}\neq 0}a_{n_{k-1},x_1} a_{n_{k-1},x_2}\dots a_{n_{k-1},x_{m_s-1}}    \nonumber \\
		\le&\frac{3\ln k}{k^{3/2}}-\frac{N}{3\cdot2^{\tau}}k^{2\tau-3}<0,  \nonumber
	\end{align}
	where the second inequality is due to (\ref{go949}) and the fourth inequality is due to (\ref{gong29}). This shows that on the set $\Omega_0\cap A$, for all large $k$, if $a_{n_{k-1},0}>2/3$, then $a_{n_{k},0}<a_{n_{k-1},0}$. Thus, by (\ref{goxin}) and (\ref{gong34}), we find that on the set $A$, except for a null set, $\liminf a_{n0}<1$, which completes the proof of Step IV.
			
\subsubsection{Proof of Step V}
	Before we prove Step V, we state
	two key lemmas, whose proofs are postponed to Section \ref{prfles}.
	\begin{lemma} \label{lemma416}
		We assume that Assumptions \ref{as1}$\sim$\ref{as2} hold. Let $j_{0}$ be a generator of the modular group of order $d$, and
		\begin{equation}  \label{gong21}
			\Omega_{K}=\bigcap_{k=K}^\infty\bigcap_{j=0}^{d-1}\left\{\max_{n_{k-1}<n\le n_k}\left|\frac1{n+n_0}\sum_{t=n_{k-1}+1}^n Q_{tj}\right| \le 2k^{-3/2}\ln k\right\}.
		\end{equation}
		If $\omega\in\Omega_{K}$ and for some large $k>K$,
		\begin{enumerate}[label=(\roman*):]
			\item  $a_{n_{k},j_0}> \mu>0$, then, for all $t\in\{0,1,\dots, d-1\}$ and $n\in [tn_k,(2d-1)n_k]$,
			$$a_{n,j_0^{*[tm_1-(t-1)]}}>p^{t-1}_1\mu^{(m_1-1)t+1}/(2d)^{m_1t+1},$$ where $j_0^{*[tm_1-(t-1)]}:=j_0\ast j_0\ast\dots\ast j_0=([tm_1-(t-1)]j_0)\mod d$.
			\item $a_{n_k,j}:=\min\{a_{n_k,j_0},a_{n_k,0}\}>\mu>0$,
			then, for all $t\in\{1,2,\dots, d\}$ and $n\in [tn_k,(2d-1)n_k]$, there is
			$$a_{n,j^{\ast t}}>p^{t-1}_1\mu^{(m_1-1)t+2-m_1}/(2d)^{m_1t+1-m_1},$$ where
			$j^{*t}:=j\ast j\ast\dots\ast j=(tj)\mod d$.
		\end{enumerate}

	\end{lemma}
	\begin{remark}\label{remark55}
	This is the only step that we need $d$ to be a prime number. When $d$ is a prime number, for all $j_0\in\{1,2,\dots, d-1\}$, $j_0$ is a generator of the modular group. Since $(m_1-1,d)=1$, $\{[(tm_1-(t-1))j_0] \mod d,\ t=0,1,\dots, d-1\}=\{0,1,\dots, d-1\}$.
	\end{remark}

	The conclusion (i) of Lemma \ref{lemma416} is obtained under the condition that the order of $j_0$ is relatively prime to $m_1-1$. This condition ensures that $\{[(tm_1-(t-1))j_0]\mod d,\ t=0,1,\dots, d-1\}=\{0,1,\dots, d-1\}$ under the random multiple-drawing urn model, which leads to the generator generating all elements in a finite number of stages by induction. Actually, by (ii) of Lemma \ref{lemma416}, it is crucial that the identity element, 0, can be generated by $j_0$. In the sequel, for any group of order $d$, we still assume $(m_1-1,d)=1$ to ensure that the order of the subgroup generated by any element of the group mutually prime to $m_1-1$, so that the identity element of the subgroup can be generated by any of the elements after sufficiently many stages.
		
	\begin{lemma}   \label{lem1} If for any $\varepsilon>0$, $\hat\mu>0$, and some large $k>K$, such that $3k^{-1/2}\ln k<\hat\mu^{m-1}\varepsilon/d$, $\min_{j}a_{n_{k},j}>\hat\mu$, and
	\begin{eqnarray}
		\max_{n_{k}<n<n_{k+1}}\left|\frac{1}{n+1}\sum_{t=n_{k}+1}^{n}Q_{tj}\right|\le2k^{-3/2}\ln k, \label{gong2140}
	\end{eqnarray}
		then we have that
	\begin{eqnarray}
		\mbox{if}\ \min_{j}a_{n_{k},j}\ge\frac{1}{d}-\varepsilon,\ \mbox{then}\ \min_{j}a_{n_{k+1},j}>\frac{1}{d}-2\varepsilon, \label{gong7}
	\end{eqnarray}
		and
	\begin{eqnarray}
		\mbox{if}\ \min_{j}a_{n_{k},j}\le\frac{1}{d}-\varepsilon,\ \mbox{then}\ \min_{j}a_{n_{k+1},j}-\min_{j}a_{n_{k},j}>\frac{\hat\mu^{m-1}\varepsilon}{k} \label{gong8}.
	\end{eqnarray}
	On the contrary,
	\begin{eqnarray}
		\mbox{if}\ \max_{j}a_{n_{k},j}\le\frac{1}{d}+\varepsilon,\ \mbox{then}\ \max_{j}a_{n_{k+1},j}<\frac{1}{d}+2\varepsilon; \label{gong9}
	\end{eqnarray}
	\begin{eqnarray}
		\mbox{if}\ \max_{j}a_{n_{k},j}\ge\frac{1}{d}+\varepsilon,\ \mbox{then}\ \max_{j}a_{n_{k+1},j}-\max_{j}a_{n_{k},j}<-\frac{\hat\mu^{m-1}\varepsilon}{k} \label{gong10}.
	\end{eqnarray}
	\end{lemma}
	Now, with the aid of Lemmas \ref{lemma416} and \ref{lem1}, we prove Step V.
	By Step IV, we know that
	$$\liminf\limits_{n\to\infty}a_{n,0}<1,\ \ a.s..$$
	Therefore,\ there exists an integer $j_{0}\in\{1,2,\dots, d-1\}$ such that
	$$\limsup\limits_{n\to\infty}a_{n,j_{0}}>0,\ \ a.s..$$
	Then, we can choose a constant $\mu>0$ and a sequence
	$\{n_{k_t}\}$ such that $a_{n_{k_t},j_0}>\mu$. Because $d$ is a prime number, $j_{0}$ is a generator of the group $\mathbb{Z}(d,+)$.
	By Lemma \ref{lemma416} and by choosing a constant $t_0$ such that $k_{t_0}>K$, we have that there exists a constant $\hat\mu>0$ such that
	\begin{equation}
		\min_j a_{n_{k_{t_0}},j}>2\hat\mu.
	\end{equation}
	We note that
	\begin{equation}  \label{goo1}
		P(\lim\limits_{K\to\infty}\Omega_{K})=1.
	\end{equation}
	Then, we find that
	\begin{equation}   \label{xin0}
		\limsup\limits_{n\to\infty}\min\limits_{j\in\{0,1,\dots, d\}}a_{nj}>\hat\mu,\ \ a.s.. 
	\end{equation}
		
For any $\varepsilon\in(0,1/2)$, by ({\ref{xin0}}), there exists $\hat\mu$ such that for infinitely many $k$, $3k^{-1/2}\ln k<\hat\mu^{m-1}\varepsilon/d$ and $\min_{j}a_{n_{k},j}\ge\hat\mu$. By (\ref{goo1}), except for a null set, for all sufficiently large $k$, (\ref{gong2140}) holds. From Lemma {\ref{lem1}}, if $\min_{j}a_{n_{k},j}<\frac{1}{d}-\varepsilon$, then $\min_{j}a_{n_{k+1},j}-\min_{j}a_{n_{k},j}\ge\hat\mu^{m-1}\varepsilon/k$. We note that for all $m$, $\sum_{k=m}^{\infty}\hat\mu^{m-1}\varepsilon/k=\infty$, there must be some $k_{0}$ such that $\min_{j}a_{n_{k_{0}},j}\ge\frac{1}{d}-\varepsilon$. Then, by Lemma {\ref{lem1}}, for all $l>k_{0}$, we have that
	$$\frac{1}{d}\ge\min_{j}a_{n_{l},j}\ge\frac{1}{d}-2\varepsilon.$$
	This shows that for all $j$,
	$$a_{n_{k},j}\to\frac{1}{d},\ \ as\ k\to \infty,\ \ \ a.s..$$
	By the fact that
	$$\max\limits_{n_{k-1}<n\le n_{k}}|a_{n,j}-a_{n_{k-1},j}|\le\frac{4}{k},$$
	we have that for all $j$,
	$$a_{n,j}\to\frac{1}{d},\ \ as\ n\to \infty,\ \ \ a.s..$$
	This completes the proof of Step V.

\subsubsection{Proof of Step VI}
	By the fundamental theorem of arithmetic, i.e., every integer greater than 1 can be represented uniquely as a product of prime numbers, the proof of Step VI can be generalized to prove the following lemma.
	\begin{lemma}\label{leem2}
		If Theorem \ref{thm1} holds for integers $d_{1}$ and $d_{2}$, then Theorem \ref{thm1} holds for $d=d_{1}d_{2}$.
	\end{lemma}
	\begin{proof}[Proof of Lemma \ref{leem2}]
	Since 1 is a generator for all $G\in(\mathbb{Z}_d,+)$, and according to Lemma \ref{lemma416}, \ref{lem1} and the proof of Step V, for the MAU model with group $G$ of order $d$, we need only show that
	\begin{equation}  \label{gong4}
			\limsup_{n\to \infty}a_{n1}>0.
	\end{equation}		
	We relabel the balls as $mod(d_{1})$; that is, we relabel the balls with numbers $\{1,\ d_{1}+1,\ \dots, \ (d_{2}-1)d_{1}+1\}$ as 1, we relabel the balls with number $\{2,\ d_{1}+2,\ \dots, \ (d_{2}-1)d_{1}+2\}$ as 2 and so on. We denote the urn composition of the relabeled model as $\mathbf{S}^{(1)}_{n}$ and for all $j\in\{0,1,\dots, d-1\}$, $a^{(1)}_{nj}=S^{(1)}_{nj}/(n+n_{0})$. Note that the urn initially contains balls labeled as the generator of $G$ of order $d$. Since there is no common factor greater than 1 between the generator and $d$, the urn still contains the balls labeled as the generator of the group with order $d_1$ after relabeling. By assumptions of Lemma \ref{leem2}, we have that
		$$\lim_{n\to\infty}a^{(1)}_{nj}=\frac{1}{d_{1}},$$
	and when $n$ is large so that
		$\min\limits_{j}a^{(1)}_{nj}\ge\frac{1}{2d_{1}},$
	for each such $n$ and $j=0,1,\dots, d_{1}-1$,
		$$a^{(1)}_{nj}=a_{nj}+a_{n,d_{1}+j}+\cdots+a_{n,(d_{2}-1)d_{1}+j}.$$
	Therefore, for each $j=0,1,\dots, d_{1}-1$, there is an $l_{j}$ such that
	\begin{equation}\label{go104}
			a_{n,l_{j}d_{1}+j}\ge\frac{1}{2d}.
	\end{equation}
	Analogously, for each large $n$ and $i=0,1,\dots, d_{2}-1$, there exists an integer $t_{i}$ such that
	\begin{equation} \label{gong1}
			a_{n,t_{i}d_{2}+i}\ge\frac{1}{2d}.
	\end{equation}
	We define the order of $\tilde{j}:=l_jd_1+j$ as $\tilde{d}$. Since $(m_1-1,d)=1$, $(m_1-1,\tilde{d})=1$. From (i) of Lemma (\ref{lemma416}) and (\ref{go104}), we have that there exists a sequence $\{n_{k_t},\ t=1,2,\dots\}$ and $\alpha>0$, such that
		$$\min\{a_{n_{k_t},l_jd_1+j},a_{n_{k_t},0}\}>\alpha>0$$
	for all $t$. By noting that $[d_{2}(l_{j}d_{1}+j)]\mod(d)=(jd_{2})\mod(d)$ and by using the method in the proof of (ii) of Lemma (\ref{lemma416}), we find that for all $t$ and $n\in[d_{2}n_{k_t},(2d_{2}-1)n_{k_t}]$,
	\begin{equation}  \label{gong2}
		a_{n,jd_{2}}\ge\frac{p^{d_2-1}_1\alpha^{(m_1-1)d_2+2-m_1}}{(2d_2)^{m_1d_2+1-m_1}},\ \ a_{n,0}>\frac{\alpha}{2d_2}.
	\end{equation}
	By ({\ref{gong1}}), for all $n\in[d_{2}n_{k_t},(2d_{2}-1)n_{k_t}]$,
	\begin{equation}
			a_{n,t_{1}d_{2}+1}\ge\frac{1}{2d}.
	\end{equation}
	We choose $j=d_{1}-t_{1}$ in ({\ref{gong2}}) for $n\in[(d_{2}+1)n_{k_t},(2d_{2}-1)n_{k_t}]$,
	\begin{equation}  \label{gong3}
		a_{n,1}\ge\frac{p^{d_2}_1\alpha^{(m_1-1)d_2}}{2d(2d_2)^{m_1d_2+2-m_1}},
	\end{equation}
	which can be obtained by the same method as (\ref{eqq99}). Hence, the assertion ({\ref{gong4}}) is proved, and thus, the proof of the lemma is completed.
	\end{proof}
		
\subsection{Proof of Theorem \ref{thm4}}
	Before proving the main theorem, we first give the following two useful lemmas, whose proofs are also postponed to Section \ref{prfles}.
	\begin{lemma}\label{lemma410}
		We assume that Assumptions \ref{as1}$\sim$\ref{as2} hold. Let $g_j$ be an element of the Abelian group of order $d$ and the order of $g_j$ be $d_1$. We define
			$$\widetilde{\Omega}_{K}=\bigcap_{k=K}^\infty\bigcap_{j=0}^{d-1}\left\{\max_{n_{k-1}<n\le n_k}\left|\frac1{n+n_0}\sum_{t=n_{k-1}+1}^n Q_{t,g_{j}}\right| \le 2k^{-3/2}\ln k\right\}.$$
		If $\omega\in\widetilde{\Omega}_{K}$ and for some large $k>K$,
		\begin{enumerate}[label=(\roman*):]
		\item  $a_{n_{k},g_{j}}> \mu>0$, then, for all $t\in\{0,1,\dots, d-1\}$ and $n\in [tn_k,(2d-1)n_k]$,
				$$a_{n,g_{j}^{*[tm_1-(t-1)]}}>p^{t-1}_1\mu^{(m_1-1)t+1}/(2d)^{m_1t+1},$$ 
		where $g_{j}^{*[tm_1-(t-1)]}:=g_j\ast g_j\ast\dots\ast g_j$.						
		\item $a_{n_k,\widetilde{g}}:=\min\{a_{n_k,g_j},a_{n_k,g_0}\}>\mu>0$, then, for all $t\in\{1,2,\dots, d\}$ and $n\in [tn_k,(2d-1)n_k]$, there is
		$$a_{n,g^{\ast t}_j}>p^{t-1}_1\mu^{(m_1-1)t+2-m_1}/(2d)^{m_1t+1-m_1},$$ 
		where $g^{*t}_j:=g_j\ast g_j\ast\dots\ast g_j$.
		\end{enumerate}
	\end{lemma}
		
	\begin{lemma} \label{lemma418}
		Let	
		$G$ be a finite Abelian group, and $H$ be a cyclic group generated by a non-identity element of $G$. If the urn initially contains balls whose labels include the generators of $G$, then, after the elements in $G$ are relabeled to $G/H$ by the left coset of $H$, the urn contains balls whose labels include the generators of $G/H$.	
	\end{lemma}
	With the aid of Lemmas \ref{lemma410} and \ref{lemma418}, we prove Theorem \ref{thm4} in the following.
		
\begin{proof}[Proof of Theorem \ref{thm4}]
		For a finite Abelian group $G=\{g_{0},g_{1},\dots, g_{d-1}\}$, we assume that $g_{0}$ is the identity element. Since the identity element of the group is unique, each element has an inverse element and the inverse element is unique, then the proof procedures from Steps I-IV in Section \ref{sec3.1} are still applicable by replacing the original zero element with the identity element of the Abelian group $G$. Similarly, Lemma \ref{lem1} also holds for Abelian group $G$. Thus, we only need to prove that
		\begin{equation} \label{gong419}
			\limsup\limits_{n\to\infty}\min\limits_{j\in\{0,1,\dots, d-1\}}a_{n,g_{j}}>0,\ a.s..
		\end{equation}
				
		From Steps I-IV in Section \ref{sec3.1}, we have that
		$$\liminf\limits_{n\to\infty}a_{n,g_{0}}<1,\ \ a.s..$$
		Thus, there exist positive constants $\nu$ and $j\in\{1,2,\dots, d-1\}$, such that
		\begin{equation} \label{gong416}
			\limsup\limits_{n\to\infty}a_{n,g_{j}}>\nu>0\ \ a.s..
		\end{equation}
		We define $H$ as a cyclic group generated by $g_{j}$; then, the order of $H$ is greater than 1, and $d_{1}:={|H|}$ is a divisor of $d$. By (\ref{gong416}), Lemma \ref{lemma410} and the fact that $P(\lim\limits_{K\to\infty}\widetilde{\Omega}_{K})=1$ (from Step IV in Section \ref{sec3.1}), we have that there exists a constant $\nu_{1}>0$, such that
		\begin{equation}
			\limsup\limits_{k\to\infty}\min_{n\in [d_{1}n_{k},(2d_1-1)n_k]}\min_{t\in\{0,1,\dots, d_{1}-1\}}a_{n,g^{\ast [tm_1-(t-1)]}_{j}}>\nu_{1}\ \ a.s..
		\end{equation}
		We note that the order of $g_j$ is $d_1$, and $(m_1-1,d)=1$. Then $(m_1-1,d_1)=1$; thus $$\{g^{\ast [tm_1-(t-1)]}_{j},\ t=0,1,\dots, d_1-1\}=\{g^{\ast t}_j,\ t=1,2,\dots, d_1\}.$$
		Hence, we have
		\begin{equation} \label{gong417}
		\limsup\limits_{k\to\infty}\min_{n\in [d_{1}n_{k},(2d_1-1)n_k]}\min_{t\in\{1,2,\dots, d_{1}\}}a_{n,g^{\ast t}_{j}}>\nu_{1}\ \ a.s..
		\end{equation}
		Since the identity element of subgroup $H$ is also the identity element of group $G$, we note that $g^{\ast d_1}_j=g_0$.	
		If $d_{1}=d$, then $G$ is the cyclic group generated by $g_{j}$, and from (\ref{gong417}), (\ref{gong419}) holds obviously. Otherwise, if $d_{1}<d$, we prove (\ref{gong419}) by induction in the following.
				
		By Corollary \ref{coro1}, we obtain the almost sure convergence of the GAU model when $d$ is a prime number. We suppose that Theorem \ref{thm4} holds if $d$ has less than $k(k\ge2)$ prime factors; we consider the case that $d$ has $k$ prime factors. Since $d_{1}<d$,
		$H$ is a nontrivial subgroup of $G$, and
		by Lagrange's theorem in Group Theory,
		there is $(d_{1},d)=d_{1}$. By the equivalence relation of elements on $G$ by $H$, for all $i,k\in\{0,1,\dots, d-1\}$,
		$$\ g_{i}\sim g_{k},\ \  \mbox{ if\ and\ only\ if\ there exists\ } h\in H, \mbox {such that}\ g_{i}=g_{k}h.$$
		Thus, each element in $G$ is equivalent to the other $d_{1}-1$ elements. Let $d_2:=d/{d_1}$. We assume that $$\widetilde{g}_{11}\sim\widetilde{g}_{12}\sim\dots\sim\widetilde{g}_{1,d_{1}},$$ $$\widetilde{g}_{21}\sim\widetilde{g}_{22}\sim\dots\sim\widetilde{g}_{2,d_{1}},$$
		$$\vdots$$
		$$\widetilde{g}_{d_{2},1}\sim\widetilde{g}_{d_{2},2}\sim\dots\sim\widetilde{g}_{d_{2},d_{1}},$$ where $\{\widetilde{g}_{i,k},\ i=1,\dots, d_{2},\ j=1,\dots, d_{1}\} $ is a relabeling of $G$. Define by  $\widetilde{G}$ the set of left cosets of $H$ with respect to $G$. Since $G$ is an Abelian group, then $\widetilde{ G}:=G/H=\{\widetilde{g}_{11}H,\widetilde{g}_{21}H,\dots, \widetilde{g}_{d_{2},1}H\}$ is a factor group, and the order of $\widetilde{G}$ is $d_2$. Since $d$ has $k$ prime factors and $1<d_1<d$, $d_2$ has at most $k-1$ prime factors.
				
		The proportion of balls labeled $\widetilde{g}_{j1}H$ after $n$ steps is denoted by $a^{(1)}_{n,\widetilde{g}_{j1}H}$, $j=1,2,\dots, d_{2}$. By Lemma \ref{lemma418} and the induction hypothesis, we have that
		$$\lim\limits_{n\to\infty}a^{(1)}_{n,\widetilde{g}_{j1}H}=\frac{1}{d_{2}}, \ \ a.s.,\ \ \mbox{for\ all }j\in\{1,2,\dots, d_{2}\}.$$
		Thus, for all large $n$, we have that
		$$\min\limits_{j}a^{(1)}_{n,\widetilde{g}_{j1}H}\ge\frac{1}{2d_{2}}\ \ a.s..$$
		For each $n$ and $j=1,2,\dots, d_{2}$,
		$$a^{(1)}_{n,\widetilde{g}_{j1}H}=a_{n,\widetilde{g}_{j1}}+a_{n,\widetilde{g}_{j2}}+\cdots+a_{n,\widetilde{g}_{j,d_{1}}}.$$
		Therefore, for each $j=1,2,\dots, d_{2}$, there exists $l_{j}\in\{1,2,\dots, d_{1}\}$ such that for all large $k$,
		\begin{equation}\label{gong07}
		\min\limits_{j}a_{n,\widetilde{g}_{j,l_{j}}}\ge\frac{1}{2d}\ \ a.s.
		\end{equation}
		holds for all $n\in[d_{1}n_{k},(2d_{1}-1)n_{k}].$
		By the equivalence between elements, for all $j\in\{1,2,\dots, d_{2}\}$, $i\in\{1,2,\dots, d_{1}\}$, there exists $h_{i,l_j}\in H$, such that $\widetilde{g}_{j,i}=\widetilde{g}_{j,l_{j}}h_{i,l_j}$. Thus, we have
		\begin{eqnarray*}
			&&\limsup_{k\to\infty}\min_{n\in[(d_{1}+1)n_{k},(2d_1-1)n_k]}\min\limits_{j,i}a_{n,\widetilde{g}_{j,i}}\\
			&\ge&	\frac{p_1}{2d_1}\limsup_{k\to\infty}\min_{n\in[d_{1}n_{k},(2d_1-1)n_k]}\min\limits_{j,i}a_{n,\widetilde{g}_{j,l_{j}}}a_{n,h_{i,l_j}} a^{m_1-2}_{n,g_0}\\
			&\ge&	\frac{p_1}{2d_1}\limsup_{k\to\infty}\min_{n\in[d_{1}n_{k},(2d_1-1)n_k]}\min\limits_{j}a_{n,\widetilde{g}_{j,l_{j}}}\min\limits_{j,i}a_{n,h_{i,l_j}}\cdot a^{m_1-2}_{n,g_0}\\
			&\ge&\frac{p_1}{2d_{1}}\limsup_{k\to\infty}\min_{n\in[d_{1}n_{k},(2d_1-1)n_k]}\min\limits_{j,i}a_{n,h_{i,l_j}}a^{m_1-2}_{n,g_0}\cdot \liminf_{k\to\infty}\min_{n\in[d_{1}n_{k},(2d_1-1)n_k]}\min\limits_{j}a_{n,\widetilde{g}_{j,l_{j}}}\\
			&\ge&\frac{p_1\nu^{m_1-1}_{1}}{4dd_1}>0,\ \ a.s.,
		\end{eqnarray*}
		where the first inequality is consistent with (\ref{eqq99}), and the third inequality can be obtained from (\ref{gong417}) and (\ref{gong07}). Thus, equation (\ref{gong419}) is proved, and hence, the proof of the theorem is complete.
	\end{proof}

\subsection{Proof of Theorem \ref{thm3}}
	We assume that $d$ is an integer with $d\ge2$ (not necessarily a prime).\ Note that for all $j\in\{0,1,\dots, d-1\}$,
	\begin{eqnarray*}
	S_{n,g_{j}}&=&S_{0,g_{j}}+\sum_{q=1}^nI(Z_q=g_{j})\\
	&=&S_{0,g_{j}}+\sum_{q=1}^{n}\sum\limits_{s=1}^{c}\sum_{x_1\ast x_2\ast\dots\ast x_{m_s}=g_j}p_s a_{q-1,x_1}a_{q-1,x_2}\dots a_{q-1,x_{m_s}}+\sum_{q=1}^n Q_{q,g_j}.
	\end{eqnarray*}
	Thus, we have that
	\begin{eqnarray*} \label{gong5}
		&&S_{n,g_j}-n/d\\
		&=&S_{0,g_j}+\sum_{q=1}^n[I(Z_q=g_j)-1/d]\\
		&=&S_{0,g_j}+\sum\limits_{q=1}^{n}\sum\limits_{s=1}^{c}\sum_{x_1\ast x_2\ast\dots\ast x_{m_s}=g_j}p_s(a_{q-1,x_1}-1/d)\cdots (a_{q-1,x_{m_s}}-1/d)+\sum_{q=1}^n Q_{q,g_j}.
	\end{eqnarray*}
	We define ${\bf Q}_{n}=\sum\limits_{q=1}^{n}(Q_{q,g_0},Q_{q,g_1},\dots, Q_{q,g_{d-1}})$. In the decomposition of ${\bf S}_{n}$,\ ${\bf Q}_{n}$ forms a martingale vector sequence. Thus, by Theorem \ref{thm4} and the martingale CLT \citep{r30}, we know that
	\begin{align}\label{lm0}
		\frac{1}{\sqrt{n}}{\bf Q}_{n} \to N_{d}(0,\Sigma)\ \ \ \mbox{in\  distribution.}
	\end{align}
	Thus, from the decomposition of ${\bf S}_n$, what remains for the proof of Theorem \ref{thm3} is to show that for all $j\in\{0,1,\dots, d-1\}$,
	\begin{equation} \label{gong6}
		\frac{1}{\sqrt{n}}\sum\limits_{q=1}^{n}\sum\limits_{s=1}^{c}\sum_{x_1\ast x_2\ast\dots\ast x_m=g_j}p_s(a_{q-1,x_1}-1/d)\cdots (a_{q-1,x_{m_s}}-1/d)\to 0,\ \ \ \mbox{in \ probability}.
	\end{equation}
	Note that by the inequality of arithmetic and geometric means and the properties of group element multiplication, we have that for all $s\in\{1,2,\dots, c\}$,
	\begin{eqnarray*}
		&&\sum_{x_1\ast x_2\ast\dots\ast x_{m_s}=g_j}(a_{q-1,x_1}-1/d)(a_{q-1,x_2}-1/d)\dots (a_{q-1,x_{m_s}}-1/d)\\
		&\le&\sum_{x_1\ast x_2\ast\dots\ast x_{m_s}=g_j}|(a_{q-1,x_1}-1/d)(a_{q-1,x_2}-1/d)\dots (a_{q-1,x_{m_s}}-1/d)|\\
		&\le&d^{{m_s}-2}\sum\limits_{t=0}^{d-1}|(a_{q-1,g_{t}}-1/d)|^{m_s}.
	\end{eqnarray*}
	Since both $d$ and $m_s$ are finite integers and $2\le m_i\le C$ for all $i=1,2,\dots$, to obtain (\ref{gong6}), we simply need to prove that
	\begin{equation} \label{gong27}
	\frac{1}{\sqrt{n}}\sum\limits_{q=1}^{n}\sum\limits_{t=0}^{d-1}|(a_{q-1,g_t}-1/d)|^{2}\to0,\ \ \ \mbox{in\ probability}.
	\end{equation}
	For any $\varepsilon>0$, we define
		$$
	H_{q,n}(\varepsilon)=\bigcap_{i=q}^{n}\left\{\max\limits_{t}\left|a_{i,g_t}-\frac{1}{d}\right|\le \varepsilon\right\}.
		$$
	By Theorem \ref{thm4},\ we know that as $q=q_{n}\to \infty$,
	$$ P(H_{q_{n},n}(\varepsilon))\to 1.$$
	We choose an integer sequence $K_{n}\to \infty$ such that $K_{n}/\sqrt{n}\to 0$.\ In the sequel,\ $H_{q_{n},n}(\varepsilon)$ is simplified to $H_{q_{n},n}.$\ Let $q_{n}=\lfloor K_{n}\rfloor$. It follows that ({\ref{gong27}}) is equivalent to
	\begin{equation} \label{eq:1}
		\frac{1}{\sqrt{n}}\sum\limits_{q=q_{n}+1}^{n}\sum\limits_{t=0}^{d-1}\mathbb{E}\left\{I(H_{q_{n},n})\left|a_{q-1,g_t}-\frac{1}{d}\right|^{2}\right\}\to 0.
	\end{equation}
	We choose $\varepsilon<1/(4d)$.\ Then, for all large $q$,\ we have that $(q+n_{0}-1+d\varepsilon)^{2}\le(q+n_{0}-1-1/2)(q+n_{0})$.\ Now,\ for $q\in[q_{n}+1,n]$,\ we have that
	\begin{eqnarray*}
		a_{q,g_t}-\frac{1}{d}&=& \frac{q+n_{0}-1}{q+n_{0}}\left(a_{q-1,g_t}-\frac{1}{d}\right)\\
		&&+\frac{1}{q+n_{0}}\left\{\sum\limits_{s=1}^{c}p_s\sum\limits_{x_1\ast\dots\ast x_{m_s}=g_t}\left(a_{q-1,x_1}-\frac{1}{d}\right)\dots\left(a_{q-1,x_{m_s}}-\frac{1}{d}\right)+Q_{q,g_t}\right\}.
	\end{eqnarray*}
	Thus,
	\begin{eqnarray*}
		&&\mathbb{E}\left\{I(H_{q_{n},q-1})\left(a_{q,g_t}-\frac{1}{d}\right)^{2}\right\}\\
		&=&\frac{(q+n_{0}-1)^{2}}{(q+n_{0})^{2}}\mathbb{E}\left\{I(H_{q_{n},q-1})\left(a_{q-1,g_t}-\frac{1}{d}\right)^{2}\right\}+\frac{2(q+n_{0}-1)}{(q+n_{0})^{2}}\mathbb{E}\left\{I(H_{q_{n},q-1})\sum\limits_{s=1}^{c}p_s\cdot\notag\right.
		\\
		\phantom{=\;\;}&&
		\left.
		\sum\limits_{x_1\ast x_2\ast\dots\ast x_{m_s}=g_t}\left(a_{q-1,g_t}-\frac{1}{d}\right)\left(a_{q-1,x_1}-\frac{1}{d}\right)\dots \left(a_{q-1,x_{m_s}}-\frac{1}{d}\right)\right\}\\
		&&+\frac{1}{(q+n_0)^{2}}\mathbb{E}\left\{I(H_{q_{n},q-1})\left[\sum\limits_{s=1}^{c}p_s\sum\limits_{x_1\ast \dots\ast x_{m_s}=g_t}\left(a_{q-1,x_1}-\frac{1}{d}\right)\dots\left(a_{q-1,x_{m_s}}-\frac{1}{d}\right)\right]^{2}\right\}\\
		&&+\frac{1}{(q+n_{0})^{2}}\mathbb{E}\left\{I(H_{q_{n},q-1})Q^{2}_{q,g_t}\right\}\\
		&\le&\frac{(q+n_{0}-1)^{2}}{(q+n_{0})^{2}}\mathbb{E}\left\{I(H_{q_{n},q-1})\left(a_{q-1,g_t}-\frac{1}{d}\right)^{2}\right\}+\frac{1}{(q+n_{0})^{2}}\\
		&&+\frac{2(q+n_{0}-1)\varepsilon}{(q+n_{0})^{2}}\mathbb{E}\left\{I(H_{q_{n},q-1})\sum\limits_{s=1}^{c}p_s\sum\limits_{x_1\ast\dots\ast x_{m_s}=g_t}\left|\left(a_{q-1,x_1}-\frac{1}{d}\right)\dots\left(a_{q-1,x_{m_s}}-\frac{1}{d}\right)\right|\right\}\\
		&&+\frac{1}{(q+n_0)^{2}}\mathbb{E}\left\{I(H_{q_{n},q-1})\left[\sum\limits_{s=1}^{c}p_s\sum\limits_{x_1\ast\dots\ast x_{m_s}=g_t}\left|\left(a_{q-1,x_1}-\frac{1}{d}\right)\dots\left(a_{q-1,x_{m_s}}-\frac{1}{d}\right)\right|\right]^{2}\right\}\\
		&\le&\frac{(q+n_{0}-1)^{2}}{(q+n_{0})^{2}}\mathbb{E}\left\{I(H_{q_{n},q-1})\left(a_{q-1,g_t}-\frac{1}{d}\right)^{2}\right\}+\frac{1}{(q+n_{0})^{2}}\\
		&&+\frac{2(q+n_{0}-1)\varepsilon}{(q+n_{0})^{2}}\mathbb{E}\left\{I(H_{q_{n},q-1})\sum\limits_{s=1}^{c}p_sd^{m_s-2}\sum\limits_{t=0}^{d-1}\left|a_{q-1,g_t}-\frac{1}{d}\right|^{m_s}\right\}\\
		&&+\frac{1}{(q+n_0)^{2}}\mathbb{E}\left\{I(H_{q_{n},q-1})\left[\sum\limits_{s=1}^{c}p_sd^{2m_s-4}\sum\limits_{t=0}^{d-1}\left|a_{q-1,g_t}-\frac{1}{d}\right|^{m_s}\right]^{2}\right\}\\
		&\le&\frac{(q+n_{0}-1)^{2}}{(q+n_{0})^{2}}\mathbb{E}\left\{I(H_{q_{n},q-1})\left(a_{q-1,g_t}-\frac{1}{d}\right)^{2}\right\}+\frac{1}{(q+n_{0})^{2}}\\
		&&+\frac{2(q+n_{0}-1)\varepsilon}{(q+n_{0})^{2}}\mathbb{E}\left\{I(H_{q_{n},q-1})\sum\limits_{s=1}^{c}p_sd^{m_s-2}\varepsilon^{m_s-2}\sum\limits_{t=0}^{d-1}\left(a_{q-1,g_t}-\frac{1}{d}\right)^{2}\right\}\\
		&&+\frac{d\varepsilon^{2}}{(q+n_0)^{2}}\sum\limits_{t=0}^{d-1}\mathbb{E}\left\{I(H_{q_{n},q-1})\sum\limits_{s=1}^{c}p_sd^{2m_s-4}\varepsilon^{2m_s-4}\left(a_{q-1,g_t}-\frac{1}{d}\right)^{2}\right\}\\
		&\le&\frac{(q+n_{0}-1)^{2}}{(q+n_{0})^{2}}\mathbb{E}\left\{I(H_{q_{n},q-1})\left(a_{q-1,g_j}-\frac{1}{d}\right)^{2}\right\}+\frac{1}{(q+n_{0})^{2}}\\
		&&+\frac{2(q+n_{0}-1)\varepsilon+d\varepsilon^{2}}{(q+n_{0})^{2}}\sum\limits_{t=0}^{d-1}\mathbb{E}\left\{I(H_{q_{n},q-1})\left(a_{q-1,g_t}-\frac{1}{d}\right)^{2}\right\},
	\end{eqnarray*}
	where the penultimate inequality is because we can choose a small enough $\varepsilon$ to satisfy $d\varepsilon<1$ for all fixed $d$. Then, we have
	\begin{eqnarray*}
		&&\sum\limits_{t=0}^{d-1}\mathbb{E}\left\{I(H_{q_{n},n})\left(a_{q,g_t}-\frac{1}{d}\right)^{2}\right\}\\
		&\le&\frac{(q+n_{0}-1+d\varepsilon)^{2}}{(q+n_{0})^{2}}\sum\limits_{t=0}^{d-1}\mathbb{E}\left\{I(H_{q_{n},q-1})\left(a_{q-1,g_t}-\frac{1}{d}\right)^{2}\right\}+\frac{d}{(q+n_{0})^{2}}\\
		&\le&\frac{q+n_{0}-1-1/2}{q+n_{0}}\sum\limits_{t=0}^{d-1}\mathbb{E}\left\{I(H_{q_{n},q-1})\left(a_{q-1,g_t}-\frac{1}{d}\right)^{2}\right\}+\frac{d}{(q+n_{0})^{2}}\\
		&\le&\left(\prod\limits_{i=q_{n}+n_{0}}^{q+n_{0}-1}\frac{i-1/2}{i+1}\right)\sum\limits_{t=0}^{d-1}\mathbb{E}\left\{I(H_{q_{n},q_{n}})\left(a_{q_{n},g_t}-\frac{1}{d}\right)^{2}\right\}+\sum\limits_{i=q_{n}+n_{0}}^{q+n_{0}-1}\frac{d}{(i+1)^{2}}\left(\prod\limits_{j=i+1}^{q+n_{0}-1}\frac{j-1/2}{j+1}\right)\\
		&\le&\frac{2(q_{n}+n_{0}-1)^{3/2}}{(q+n_{0}-1)^{3/2}}d\varepsilon^{2}+\sum\limits_{i=q_{n}+n_{0}}^{q+n_{0}-1}\frac{2di^{3/2}}{(i+1)^{2}(q+n_{0}-1)^{3/2}},
	\end{eqnarray*}
	where we define $\prod\limits_{j=q+n_{0}}^{q+n_{0}-1}\frac{j-1/2}{j+1}=1$. \ Consequently, as $n\to\infty$,
	\begin{eqnarray*}
		&&n^{-1/2}\sum\limits_{q=q_{n}+1}^{n}\sum\limits_{t=0}^{d-1}\mathbb{E}\left\{I(H_{q_{n},n})\left(a_{q,g_t}-\frac{1}{d}\right)^{2}\right\}\\
		&\le&n^{-1/2}\sum\limits_{q=q_{n}+1}^{n}\frac{2(q_{n}+n_{0}-1)^{3/2}}{(q+n_{0}-1)^{3/2}}d\varepsilon^{2}+n^{-1/2}\sum\limits_{q=q_{n}+1}^{n}\sum\limits_{i=q_{n}+n_{0}}^{q+n_{0}-1}\frac{2di^{3/2}}{(i+1)^{2}(q+n_{0}-1)^{3/2}}\\
		&\le&4dn^{-1/2}(q_{n}+n_{0}-1)\varepsilon^{2}+3dn^{-1/2}\ln n \ \to\ 0.
	\end{eqnarray*}
	Then, the proof of Theorem \ref{thm3} is complete.
		
\subsection{Proofs of Lemmas \ref{le1}, \ref{lemma416}, \ref{lem1} and \ref{lemma418}}\label{prfles}
	Next, we give the proofs of Lemmas \ref{le1} \ref{lemma416}, \ref{lem1} and \ref{lemma418}. Note that as the proof of Lemma \ref{lemma410} is consistent with that of Lemma \ref{lemma416}, the details are omitted in this section.
	\begin{proof}[Proof of Lemma \ref{le1}]
		\begin{eqnarray*}
			P\left(T_1=\infty\right)&=&P\left(\lim\limits_{n\to\infty}\bigcap\limits_{i=1}^{n}\{M_i\neq m_1\}\right)=\lim\limits_{n\to\infty}P\left(\bigcap\limits_{i=1}^{n}\{M_i\neq m_1\}\right)\\
			&=&\lim\limits_{n\to\infty}\prod\limits_{i=1}^{n}P(M_i\neq m_1)=\lim\limits_{n\to\infty}(1-p_1)^{n}=0,
		\end{eqnarray*}
		which shows that $T_1<\infty$ almost surely. Then,
		\begin{eqnarray*}
			P(T_2<\infty)&=&P(\{T_2<\infty\}\cap\{T_1<\infty\})+P(\{T_2<\infty\}\cap\{T_1=\infty\})\\
			&=&P(\{T_2<\infty\}\cap\{T_1<\infty\})=1-P(\{T_2=\infty\}\cap\{T_1<\infty\})\\
			&\ge&1-P(\{T_2>T_1\}\cap\{T_1<\infty\})=1-P(\{T_2>T_1\}\cap\{T_2\ge T_1\}|\{T_1<\infty\})\\
			&=&1-P(\{T_2>T_1\}|\{T_2\ge T_1\}\cap\{T_1<\infty\})P(T_2\ge T_1|T_1<\infty)=1
		\end{eqnarray*}	
		Thus, Lemma \ref{le1} is proved.
	\end{proof}

	\begin{proof}[Proof of Lemma \ref{lemma416}]
		We first give the proof of (i) of Lemma \ref{lemma416}. When $t=0$, for all $n\in[n_k,(2d-1)n_k]$, we have
		\begin{align}\label{gong55}
			a_{n,j_0}\ge\frac{S_{0,j_0}+\sum\limits_{i=1}^{n_{k}}I(Z_{i}=j_0)}{n+n_{0}}\ge\frac{(n_{k}+n_{0})a_{n_{k},j_0}}{n+n_{0}}\ge\frac{n_{k}+n_{0}}{(2d-1)n_{k}+n_{0}}a_{n_{k},j_0}>\frac{\mu}{2d}>\frac{\mu^{m_1}}{(2d)^{m_1+1}}.
		\end{align}
		By induction,\ the conclusion for $t\ge1$ follows from
		\begin{eqnarray}
			&&a_{n,j^{\ast [tm_1-(t-1)]}_0}\nonumber\\
			&\ge&\frac{1}{(2d-1)n_k+n_0}\sum_{i=(t-1)n_k+1}^{tn_k}I(Z_i=j^{\ast [tm_1-(t-1)]}_0)\nonumber\\
			&=&\frac{1}{(2d-1)n_k+n_0}\sum_{i=(t-1)n_k+1}^{tn_k}\left[P\left(Z_i=j^{\ast [tm_1-(t-1)]}_0|{\cal F}_{i-1}\right)+Q_{i,j^{\ast [tm_1-(t-1)]}_0}\right]\nonumber\\
			&=&\frac{1}{(2d-1)n_k+n_0}\sum_{i=(t-1)n_k+1}^{tn_k}\left[\sum\limits_{s=1}^{c}P\left(Z_i=j^{\ast [tm_1-(t-1)]}_0|M_i=m_s,{\cal F}_{i-1}\right)P\left(M_i=m_s|{\cal F}_{i-1}\right)
			\notag\right.
			\\
			\phantom{=\;\;}&&
			\left.
			+Q_{i,j^{\ast [tm_1-(t-1)]}_0}\right]\nonumber\\
			&=&\frac{1}{(2d-1)n_k+n_0}\sum_{i=(t-1)n_k+1}^{tn_k}\sum\limits_{s=1}^{c}p_s\sum\limits_{x_1+\dots+ x_{m_s}=[tm_1-(t-1)]j_0}a_{i-1,x_1}\dots a_{i-1,x_{m_s}}+O\left(\frac{\log k}{k^{1/2}}\right)\nonumber\\
			&\ge&\frac{p_1}{(2d-1)n_k+n_0}\sum_{i=(t-1)n_k+1}^{tn_k}\sum\limits_{x_1+\dots+ x_{m_1}=[tm_1-(t-1)]j_0}a_{i-1,x_1}\dots a_{i-1,x_{m_1}}+O\left(\frac{\log k}{k^{1/2}}\right)\nonumber\\
			&\ge&\frac{p_1}{(2d-1)n_k+n_0}\sum_{i=(t-1)n_k+1}^{tn_k}\sum\limits_{x_1+\dots+ x_{m_1}=[tm_1-(t-1)]j_0}a_{i-1,j_0}\dots a_{i-1,j_0}a_{i-1,[(t-1)m_1-(t-2)]j_0}\nonumber\\
			&\ge&\frac{p^{t-1}_1}{2d}\frac{\mu^{m_1-1}}{(2d)^{m_1-1}}\frac{\mu^{(m_1-1)(t-1)+1}}{(2d)^{m_1(t-1)+1}}=\frac{p^{t-1}_1\mu^{(m_1-1)t+1}}{(2d)^{m_1t+1}}\nonumber.
		\end{eqnarray}
		which completes the proof of (i) of Lemma \ref{lemma416}. We then prove (ii). (\ref{gong55}) still holds. That is, for all $n\in\{n_k,(2d-1)n_k\}$, $a_{n,j}>\mu/(2d)$. By induction, for $t\ge2$, we have
			\begin{align}
		 	&a_{n,j^{\ast t}}\nonumber\\
		 	\ge&\frac{1}{(2d-1)n_k+n_0}\sum_{i=(t-1)n_k+1}^{tn_k}I(Z_i=j^{\ast t})\nonumber\\
		 	=&\frac{1}{(2d-1)n_k+n_0}\sum_{i=(t-1)n_k+1}^{tn_k}\left[P\left(Z_i=j^{\ast t}|{\cal F}_{i-1}\right)+Q_{i,j^{\ast t}}\right]\nonumber\\
		 	=&\frac{1}{(2d-1)n_k+n_0}\sum_{i=(t-1)n_k+1}^{tn_k}\left[\sum\limits_{s=1}^{c}P\left(Z_i=j^{\ast t}|M_i=m_s,{\cal F}_{i-1}\right)P\left(M_i=m_s|{\cal F}_{i-1}\right))+Q_{i,j^{\ast t}}\right]\nonumber\\
		 	=&\frac{1}{(2d-1)n_k+n_0}\sum_{i=(t-1)n_k+1}^{tn_k}\sum\limits_{s=1}^{c}p_s\sum\limits_{x_1+\dots+ x_{m_s}=tj}a_{i-1,x_1}\dots a_{i-1,x_{m_s}}+O\left(\frac{\log k}{k^{1/2}}\right)\nonumber\\
		 	=&\frac{1}{(2d-1)n_k+n_0}\sum_{i=(t-1)n_k+1}^{tn_k}\sum\limits_{s=1}^{c}p_s\sum\limits_{x_1+ \dots+ x_{m_s}=tj}a_{i-1,x_1}a_{i-1,x_2}\dots a_{i-1,x_{m_s}}+O\left(\frac{\log k}{k^{1/2}}\right)\nonumber\\
		 	\ge&\frac{p_1}{(2d-1)n_k+n_0}\sum_{i=(t-1)n_k+1}^{tn_k}\sum\limits_{x_1+ \dots+ x_{m_1}=tj}a_{i-1,j}a_{i-1,(t-1)j}a_{i-1,0}\dots a_{i-1,0}+O\left(\frac{\log k}{k^{1/2}}\right)\nonumber\\
		 	\ge&\frac{p_1}{2d}a_{(t-1)n_k,j}a_{(t-1)n_k,(t-1)j} a_{(t-1)n_k,0}\dots a_{(t-1)n_k,0}a_{(t-1)n_k,0}\label{eqq99}\\
		 	>&\frac{p^{t-1}_1}{2d}\frac{\mu^{m_1-1}}{(2d)^{m_1-1}}\frac{\mu^{(m_1-1)(t-1)+2-m_1}}{(2d)^{m_1(t-1)+1-m_1}}=\frac{p^{t-1}_1\mu^{(m_1-1)t+2-m_1}}{(2d)^{m_1t+1-m_1}}\nonumber.
		 \end{align}
	\end{proof}
	\begin{proof}[Proof of Lemma \ref{lem1}]
		As before (by (\ref{oo1})), we know that
		\begin{align}
			&a_{n_{k+1},j}-a_{n_{k},j}\nonumber\\
			=&\frac{1}{n_{k+1}+n_{0}}\left[(n_{k}+n_{0})a_{n_k,j}+\sum_{m=n_{k}+1}^{n_{k+1}}I(Z_m=j)\right]-a_{n_{k},j}\nonumber\\
			=&\frac{1}{n_{k+1}+n_{0}}\sum_{m=n_{k}+1}^{n_{k+1}}[I(Z_m=j)-a_{n_{k},j}]\nonumber\\ 
			=&\frac{1}{n_{k+1}+n_{0}}\sum_{m=n_{k}+1}^{n_{k+1}}[P(Z_{m}=j|{\cal F}_{m-1})-a_{n_{k},j}]+O_1\left(\ln k/k^{3/2}\right)\nonumber\\
			=&\frac{n_{k+1}-n_{k}}{n_{k+1}+n_{0}}\left[\sum\limits_{s=1}^{c}p_sP(Z_{n_{k}+1}=j|M_{n_{k}+1}=m_s,{\cal F}_{n_{k}})-a_{n_{k},j}\right]+O_2\left(\ln k/k^{3/2}\right)\nonumber\\
			=&\frac{n_{k+1}-n_{k}}{n_{k+1}+n_{0}}\sum\limits_{s=1}^{c}p_s\left[P(Z_{n_{k}+1}=j|M_{n_k+1}=m_s,{\cal F}_{n_{k}})-a_{n_{k},j}\right]+O_2\left(\ln k/k^{3/2}\right), \label{gong11}
		\end{align}
		where the last step by (\ref{oo2}), $\left|O_1\left(\ln k/k^{3/2}\right)\right|\le2\ln k/k^{3/2}$, and $\left|O_2\left(\ln k/k^{3/2}\right)\right|\le3\ln k/k^{3/2}$. Now, we use ({\ref{gong11}}) to prove the Lemma. We assume that $a_{n_{k},I}=\min_{j}a_{n_{k},j}\le\frac{1}{d}-\varepsilon$ and $a_{n_{k},J}=\max_{j}a_{n_{k},j}$; then, by the fact that $a_{n_{k},J}> \frac{1}{d}+\frac{\varepsilon}{d}$,\ for any $j$,				

		\begin{eqnarray*}
			&&P\left(Z_{n_{k}+1}=j|{\cal F}_{n_k}\right)-a_{n_k,I}\\
			&=&\sum\limits_{x_1\ast x_2\ast\dots\ast x_m=j}a_{n_k,x_1}a_{n_k,x_2}\dots a_{n_k,x_m}-a_{n_k,I}\\
			&=&\sum\limits_{x_1\ast x_2\ast\dots\ast x_m=j}a_{n_k,x_1}a_{n_k,x_2}\dots a_{n_k,x_{m-1}}(a_{n_k,x_m}-a_{n_k,I})\\
			&\ge&(a_{n_k,J}-a_{n_k,I})a_{n_k,\hat x_1}a_{n_k,\hat x_2}\dots a_{n_k,\hat x_{m-1}}\\
			&\ge&(1+1/d)\hat\mu^{m-1}\varepsilon.
		\end{eqnarray*}
		The second equal sign comes from (\ref{go949}). This implies that if $k$ satisfies $3k^{-1/2}\ln k<\hat\mu\varepsilon/d$, by ({\ref{gong11}}), we know that
				$$
	\min_{j}a_{n_{k+1},j}-a_{n_{k},I}>\frac{n_{k+1}-n_{k}}{n_{k+1}+n_0}(1+1/d)\hat\mu\varepsilon-3k^{-3/2}\ln k>\hat\mu\varepsilon/k.
				$$
		Then, ({\ref{gong8}}) holds. The conclusion ({\ref{gong7}}) follows from the fact that
		\begin{eqnarray*}
			&&P\left(Z_{n_{k}+1}=j|{\cal F}_{n_k}\right)-a_{n_k,I}\\
			&=&\sum\limits_{x_{1}\ast x_{2}\ast\dots\ast x_m=j}a_{n_k,x_1}a_{n_k,x_2}\dots a_{n_k,x_m}-a_{n_k,I}\\
			&=&\sum\limits_{x_{1}\ast x_{2}\ast\dots\ast x_m=j}a_{n_k,x_1}a_{n_k,x_2}\dots a_{n_k,x_{m-1}}(a_{n_k,x_m}-a_{n_k,I})\ge0\\
		\end{eqnarray*}
		which implies
		$$\min_{j}a_{n_{k+1},j}\ge a_{n_{k},I}-3k^{-3/2}\ln k\ge \frac{1}{d}-2\varepsilon$$
		by (\ref{gong11}). Then, ({\ref{gong9}}) and ({\ref{gong10}}) follow by symmetry. The proof of the lemma is complete.
	\end{proof}
	\begin{proof}[Proof of Lemma \ref{lemma418}]
		Since any finite Abelian group can be expressed as the direct product of cyclic groups, we assume that the base of the group $G$ is $\{a_{1},a_{2},\dots, a_{s}\}$; then, it follows that
		$$G=\langle a_{1}\rangle\otimes\langle a_{2}\rangle\otimes\dots\otimes\langle a_{s}\rangle,$$
		where $\langle a_{i}\rangle$ denotes the cyclic group generated by $a_{i}$. We assume that the order of $a_{i}$ is $n_{i},\ i=1,2,\dots, s$. Then, for all $g\in G$, $g$ can be uniquely represented as
		$$g=a^{k_{1}}_{1}a^{k_{2}}_{2}\dots a^{k_{s}}_{s},\ 0\le k_{i}<n_{i},\ i=1,2,\dots, s.$$
		Assuming $H=\langle a^{k}_{s}\rangle$, we prove the conclusion of Lemma \ref{lemma418} in two cases:
		\begin{longlist}
			\item If $\langle k,n_{s}\rangle=1$, then $H=\langle a_{s}\rangle$,
			$$G/H=\{a^{k_{1}}_{1}a^{k_{2}}_{2}\dots a^{k_{s-1}}_{s-1}H,\ 0\le k_{i}< n_{i},\ i=1,2,\dots, s-1\},$$
			and $\{a_{1}H,a_{2}H,\dots, a_{s-1}H\}$ is a base of $G/H$. Since the urn initially contains balls labeled $a_{1},a_{2},\dots, a_{s}$, it is equivalent to the urn containing balls labeled $a_{1}H,a_{2}H,\dots, a_{s-1}H$ after relabeling.				
			\item If $\langle k,n_{s}\rangle=q>1$, then
			$$G/H=\{a^{k_{1}}_{1}a^{k_{2}}_{2}\dots a^{k_{s-1}}_{s-1}a^{k_{s}}_{s}H,\ 0\le k_{i}< n_{i},\ i=1,2,\dots, s-1,\ 0\le k_{s}<q\},$$
			and $\{a_{1}H,a_{2}H,\dots, a_{s}H\}$ is a base of $G/H$. Similarly, the urn contains balls labeled $a_{1}H,a_{2}H,\dots, a_{s}H$ after relabeling.
		\end{longlist}
		Then, the proof of this lemma is complete.
	\end{proof}

\section{Conclusion and possible extensions}
	In this paper. we discuss the first-order and second-order asymptotic properties of the GAU urn model proposed by P. Diaconis. For the GAU model, the growth rate of the urn composition is obtained. The conjecture of P. Diaconis that the normalized urn composition converges almost surely to a uniform distribution over group $G$ is confirmed. That is, for all $g\in G$, the limiting fraction of balls labeled by $g$ converges to $1/|G|$ almost surely. In addition, the CLT for the urn composition is shown, from which we found that the second-order asymptotic property of the urn composition is consistent with $\sqrt{n}\left(1/n\sum_{i=1}^{n}\textbf{W}_{i}-1/d{\bf1}_{d}\right)$, where $\{\textbf{W}_{i},\ i=1,2,\dots\}$ are independent and identically distributed random variables and $\mathbf{W}_{i}\sim$ \rm{Multinomial}\ $(1;1/d, \dots,  1/d)$ for all $i$.
		
	In the future, some interesting variations of the urn model proposed by Diaconis should be considered. For example, in order to relax the relatively prime condition, that is, Assumption \ref{as2}, one may consider the case where the urn initially contains more types of balls, or generalize the results of this paper to arbitrary finite groups, not just Abelian groups.

\end{document}